\documentclass[final]{article}

\usepackage[affil-it]{authblk}

\usepackage{amsmath,graphicx,latexsym,amssymb, amscd, psfrag,verbatim, amsfonts,amsthm}
\usepackage{amsfonts,bm,mathrsfs,subfigure}
\usepackage{geometry,graphicx}

\usepackage{setspace}
\usepackage{lipsum}
\usepackage{amsfonts}
\usepackage{graphicx}
\usepackage{epstopdf}
\usepackage{comment}
\usepackage{algorithm,algpseudocode}
\usepackage{hyperref}
\usepackage{booktabs}

\newtheorem{theorem}{Theorem}
\newtheorem{lemma}{Lemma}
\newtheorem{example}{Example}

\newcommand{\lm}{l_{\text{max}}}
\newcommand{\cA}{{\mathcal A}}
\newcommand{\cH}{{\mathcal H}}
\newcommand{\cK}{{\mathcal K}}
\newcommand{\cL}{{\mathcal L}}
\newcommand{\cU}{{\mathcal U}}
\newcommand{\cE}{{\mathcal E}}
\newcommand{\cM}{{\mathcal M}}

\begin{document}


\title{Sparse Hierarchical Solvers with Guaranteed Convergence}
\date{}
\author{Kai Yang\thanks{Department of Mechanical Engineering, Stanford University, Stanford, CA 94305, USA. \texttt{yangkai.stanford@gmail.com}}, 
Hadi Pouransari\thanks{Department of Mechanical Engineering, Stanford University, Stanford, CA 94305, USA. \texttt{hadip@stanford.edu}}, 
Eric Darve \thanks{Institute for Computational and Mathematical Engineering and
Department of Mechanical Engineering, Stanford University, Stanford, CA 94305, USA. \texttt{darve@stanford.edu}}
}
\maketitle
\begin{abstract}
Solving sparse linear systems from discretized PDEs is challenging. Direct solvers have in many cases quadratic complexity (depending on geometry), while iterative solvers require problem dependent preconditioners to be robust and efficient.  Approximate factorization preconditioners, such as incomplete LU factorization, provide cheap approximations to the system matrix. However, even a highly accurate preconditioner may have deteriorating performance when the condition number of the system matrix increases. By increasing the accuracy on low-frequency errors, we propose a novel hierarchical solver with improved robustness with respect to the condition number of the linear system. This solver retains the linear computational cost and memory footprint of the original algorithm.
\end{abstract}



\section {Introduction}
Solving a large-scale sparse linear system
\begin{equation}
\label{eq:Axb}
Ax=b
\end{equation}
with $A\in\mathbb{R}^{n\times n}$ is one of the most challenging tasks in scientific computing.  A typical feature of the matrix $A$ from discretized partial differential equations (PDEs) is sparsity, which implies that the number of nonzeros per row and per column is small compared with the number of unknowns $n$.  While it is easy to store the linear system in $\mathcal{O}(n)$ memory and multiply $A$ and a vector with $\mathcal{O}(n)$ computational cost, solving $Ax=b$ in $\mathcal{O}(n)$ (or close to $\mathcal{O}(n)$) memory and CPU time remains a challenging issue.

A straightforward Gaussian elimination results in $\mathcal{O}(n^2)$ storage and $\mathcal{O}(n^3)$ computational cost due to the fill-ins introduced. Improved ordering can reduce the fill-ins. For example, nested dissection ordering only takes storage $\mathcal{O}(n\log(n))$ in 2D and $\mathcal{O}(n^{4/3})$ in 3D for simple geometries and discretizations. However, the computational complexity of nested dissection ordering is $\mathcal{O}(n^{3/2})$ in 2D and $\mathcal{O}(n^2)$ in 3D, for most cases.

Instead of directly inverting the sparse matrix, iterative methods seek an approximate solution to (\ref{eq:Axb}) within a small number of iterations. Although iterative methods usually require only $\mathcal{O}(n)$ storage, the convergence is highly problem-dependent.  Preconditioners are needed in order to achieve fast convergence.  There is an extensive literature on developing preconditioners for various applications. Most preconditioners are problem specific and depend on detailed information, such as the geometry and discretization.  Although in general more information on the problem results in a more effective solver, in this paper we focus on robust and purely algebraic preconditioners.

Multigrid method \cite{Hackbusch.W2013a,Bramble.J1993a,Xu.J1992a,Briggs.W;Henson.V;McCormick.S2000a,Trottenberg.U;Oosterlee.C;Schuller.A2000a} is considered as one of the most effective solvers for discretized elliptic PDEs.  Based on the discretization of the PDEs, multigrid method efficiently eliminates the error in a multilevel fashion.   As an algebraic variant, algebraic multigrid (AMG) \cite{Brandt.A;McCoruick.S;Huge.J1985a,Brandt.A1986a,Ruge.J;Stuben.K1987a,Stuben.K2001a,Napov.A;Notay.Y2012a,Notay.Y2010a} does not require information of the discretization and constructs multilevel structures based on the sparsity pattern of the matrix.  AMG has been shown to be efficient for a wide range of applications, especially when the problem is close to elliptic PDEs.  In some classical model problems, like the Poisson equation with constant coefficients, AMG achieves linear complexity in memory and time in a way similar to multigrid method.  For example, in \cite{Napov.A;Notay.Y2012a} aggregation-based AMG is shown to be optimal for model Poisson problem under certain assumptions. 

Another purely algebraic preconditioner is the incomplete LU (ILU) preconditioner \cite{Saad.Y2003a}.  ILU resembles Gaussian elimination, except that it drops some of the fill-ins to keep the sparsity.  ILU has demonstrated its effectiveness in many applications such as convection-diffusion equations.  Although not optimal, ILU improves the performance over that of an unpreconditioned iterative solver.  Meanwhile, choosing the proper dropping rules for an efficient ILU factorization is a problem-dependent trial and error process.

Low-rank compression is another type of approaches to construct an efficient and accurate preconditioner.  This approach was used in fast multipole method (FMM) \cite{Greengard.L;Rokhlin.V1987a,Darve.E2000a,Nishimura.N2002a,Ying.L;Biros.G;Zorin.D2004a,Fong.W;Darve.E2009a,Pouransari.H;Darve.E2015a} to accelerate matrix-vector multiplication and then generalized to hierarchical matrices ($\mathcal{H}$-matrices) \cite{Hackbusch.W;Borm.S2002a,Borm.S;Grasedyck.L;Hackbusch.W2003a,Bebendorf.M2008a} for solving dense linear systems from integral equations.  There are different types of $\mathcal{H}$-matrices depending on the low-rank admissibilities and whether the low-rank bases are nested.  The hierarchical off-diagonal low-rank (HODLR) matrix \cite{Ambikasaran.S;Darve.E2013a,Aminfar.A;Ambikasaran.S;Darve.E2016a,Aminfar.A;Darve.E2016a} assumes the off-diagonal blocks to be low-rank. If the low-rank bases of these off-diagonal blocks can in addition be nested, then this matrix is called hierarchically semi-separable (HSS) \cite{Chandrasekaran.S;Gu.M;Pals.T2006a, Xia.J;Chandrasekaran.S;Gu.M;Li.X2010a,Ambikasaran.S;Darve.E2013a}.  More general than HSS, $\mathcal{H}^2$ matrices are assumed to have low-rank structure only for well-separated interactions \cite{Borm.S;Grasedyck.L;Hackbusch.W2003a}, while maintaining nested low-rank bases.  In addition to accelerating matrix-vector multiplication, these low-rank techniques have been also applied to approximately factorize the dense matrix and provide direct solvers \cite{Bebendorf.M2008a,Ambikasaran.S;Darve.E2014a,Coulier.P;Pouransari.H;Darve.E2015a}.

Although mostly applied to solving dense linear systems \cite{Greengard.L;Gueyffier.D;Martinsson.P;Rokhlin.V2009a,Kong.W;Bremer.J;Rokhlin.V2011a,Gillman.A;Young.P;Martinsson.P2012a,Corona.E;Martinsson.P;Zorin.D2015a,Coulier.P;Pouransari.H;Darve.E2015a}, hierarchical matrix techniques have also been applied to solving sparse linear systems, especially those from the discretization of second order elliptic PDEs.  For example, the hierarchical LU and Cholesky factorization developed by Hackbusch et al. \cite{Hackbusch.W1999a,Hackbusch.W;Borm.S2002a,Borm.S;Grasedyck.L;Hackbusch.W2003a,Bebendorf.M2007a} follow a tree decomposition, recursively factorizing block matrices and forming Schur complements.  More recently, hierarchical interpolative factorization \cite{Ho.K;Ying.L2015a,Ho.K;Ying.L2015b} was proposed to solve linear systems from the discretization of integral and differential equations based on elliptic operators.  Many other approaches based on low-rank properties have also been proposed to solving sparse linear systems \cite{Li.R;Saad.Y2013a,Napov.A;Li.X2016a,Sushnikova.D;Oseledets.I2016a}.

More recently, \cite{Pouransari.H;Coulier.P;Darve.E2015a} extended the approaches in \cite{Ambikasaran.S;Darve.E2014a,Coulier.P;Pouransari.H;Darve.E2015a} to solving sparse linear systems, assuming that the fill-ins during a block Gaussian elimination have low numerical rank for well-separated interactions. An approximate factorization has been introduced via extended sparsification \cite{Pouransari.H;Coulier.P;Darve.E2015a, Chandrasekaran.S;Dewilde.P;Gu.M;Lyons.W;Pals.T2006a}, which has time and memory complexity of $\mathcal{O}(n\log^21/\epsilon)$ and $\mathcal{O}(n\log1/\epsilon)$, respectively.  Here $\epsilon$ is the accuracy of the low-rank approximation.  This application results in an accurate preconditioner for the original linear system. 
Despite the high accuracy of approximate factorizations, the robustness of this type of approaches is subject to the condition number of the original problem.  In \cite{Pouransari.H;Coulier.P;Darve.E2015a}, rather than solving (\ref{eq:Axb}) directly, the following equation is solved instead:
\begin{equation}
\label{eq:approx}
 A_{\cH}x_{\cH}=b
 \qquad \text{with} \qquad
\|A_{\cH}-A\|\leq \epsilon \|A\|
\end{equation} 
Then the accuracy of the solution is subject to both $\epsilon$ and $\kappa(A)=\|A\|\|A^{-1}\|$ 
\[ \frac{\|x-x_{\cH}\|}{\|x\|}\leq C\epsilon \kappa(A) \]
where $C$ is a constant.  As the upper bound for the relative accuracy is proportional to the matrix condition number, the accuracy deteriorates when the condition number of matrix grows. This can break the linear complexity as the condition number increases with the number of unknowns in many applications.

In recent work \cite{Bebendorf.M;Bollhofer.M;Bratsch.M2013a,Bebendorf.M;Bollhofer.M;Bratsch.M2016a},  new techniques that alleviate this sensitivity are proposed.   The approximate factorization is required to preserve additional vectors, called side constraints, in order to improve the robustness of the solver.  It is proved in \cite{Bebendorf.M;Bollhofer.M;Bratsch.M2013a} that  the approximate factorization that additionally preserves low-frequency eigenvectors is spectrally equivalent to the system matrix.  \cite{Bebendorf.M;Bollhofer.M;Bratsch.M2016a} also provides adaptive choices of compression parameters for the preconditioner to be spectrally equivalent with constant side constraints.  In this paper, we follow similar ideas and propose an improved version of the algorithm in \cite{Pouransari.H;Coulier.P;Darve.E2015a}. We note that although the terminology ``hierarchical solver'' is used both in \cite{Pouransari.H;Coulier.P;Darve.E2015a} and \cite{Bebendorf.M;Bollhofer.M;Bratsch.M2013a,Bebendorf.M;Bollhofer.M;Bratsch.M2016a}, these methods are distinct. The distinction was clarified in \cite{Coulier.P;Pouransari.H;Darve.E2015a} for example (see Appendix C at the end of this reference). Essentially, \cite{Bebendorf.M;Bollhofer.M;Bratsch.M2013a,Bebendorf.M;Bollhofer.M;Bratsch.M2016a} is based on a block LU factorization, while \cite{Pouransari.H;Coulier.P;Darve.E2015a} is based on a multilevel algorithm similar to AMG.

In this paper, we first formulate the algorithm LoRaSp from \cite{Pouransari.H;Coulier.P;Darve.E2015a} in terms of block matrix factorization and analyze the error of the approximate factorization. Then, we present an estimate of the condition number of the preconditioned linear system, which implies that the performance of the hierarchical solver is restricted by its convergence on low-frequency errors. Based on this observation, we enrich the low-rank bases in the matrix compression process to improve the convergence of our hierarchical solver on low-frequency eigenvectors. We call this new method GC-eigenvector (Guaranteed Convergence based on preserving eigenvectors). Numerical examples illustrate that this modification greatly improves the robustness of the solver. In particular, the number of iterations in GMRES becomes nearly independent of the problem size, whereas the condition number increases like $O(1/h^2)$ in the benchmark problems.

GC-eigenvector requires computing some of the small eigenvectors of $A$, which is computationally expensive. To make this algorithm more practical, the constant vector, which is a computationally cheap alternative to small eigenvectors\footnote{For brevity, we will use the expression ``small'' eigenvector, to mean the eigenvector associated with a small eigenvalue.}, is instead preserved in the hierarchical solver, and investigated in numerical benchmarks in Section \ref{sec:numerics}. This method is called GC-constant. We show that GC-constant is both computationally efficient as it requires minimal additional cost compared to our reference LoRaSp hierarchical solver, and leads to the same the convergence rate for GMRES as GC-eigenvector. 

In summary, the main contributions of this paper are: 
\begin{enumerate}
\item Numerical analysis of the LoRaSp \cite{Pouransari.H;Coulier.P;Darve.E2015a} h-solver. Theorems related to the condition number of the preconditioned system are proved.
\item A modification of LoRaSp based on the smallest eigenvectors of $A$ (``low-frequency'' eigenvectors), GC-eigenvector, is proposed and benchmarked. This leads a nearly constant number of iterations in GMRES, irrespective of the condition number of $A$.
\item A computationally more efficient method called GC-constant is proposed. It leads to a preconditioner which is as efficient as GC-eigenvector (in terms of reducing the number of iterations in GMRES), while having a computational cost comparable to the original LoRaSp. No eigenvector of $A$ needs to be estimated in this method.
\item Numerical benchmarks for various elliptic PDEs are used to validate our analysis and the new algorithms.
\end{enumerate}

The remainder of this paper is organized as follows. In Section \ref{sec:sparse} and \ref{sec:hsolver}, the factorization phase and solve phase of the hierarchical solver are explained in detail using block matrix notation. In Section \ref{sec:error}, we analyze the error of the hierarchical solver due to the low-rank approximations. In Section \ref{sec:convergence}, we analyze the convergence of hierarchical solver when used a preconditioner. Based on these analysis, we propose in Section \ref{sec:improved} a new technique to improve the convergence of the hierarchical solver on low-frequency errors.  In Section \ref{sec:numerics}, we benchmark the improvement due to the new technique compared with the original hierarchical solver. We verified on these problems that the convergence becomes nearly independent of the problem size.  In the appendix we propose some approaches to let the hierarchical solvers approximately preserve vectors of interest.


%

\section{Sparse linear systems}
\label{sec:sparse}
In this section, we introduce basic concepts about sparse linear systems. Throughout the rest of this paper, we assume that $A$ is a symmetric positive definite (SPD) matrix.  Then, the condition number of $A$ is defined by the ratio of extreme eigenvalues: $\kappa(A)=\lambda_\text{max}(A)/\lambda_\text{min}(A)$. 

\subsection{Graph representation of a sparse matrix}
As the matrix $A=(a_{ij})$ is assumed to be symmetric, it corresponds to an undirected graph $G=(V,E)$, where
\[ V=\{1,2,\ldots,n\},\quad E=\{(i,j)| a_{ij}\neq0 \} \]


To better leverage the fast block matrix operations of BLAS3 \cite{Anderson.E;Bai.Z;Bischof.C;Blackford.S;Demmel.J;Dongarra.J;Du-Croz.J;Greenbaum.A;Hammerling.S;McKenney.A;others1999a}, we introduce block partitioning of matrices.  Given a partitioning $\{I_i\}_{i=1}^N$ of $V$ s.t. 
\[\cup_iI_i=V,\quad I_i\cap I_j=\emptyset, \text{ if }i\neq j\]
Let $A_{I_iI_j}$ denote a block matrix of $A$ with row indices in $I_i$ and column indices in $I_j$.  A new graph with less vertices, denoted by $G_b$, is then defined as $G_b=(V_b,E_b)$ with
\[V_b=\{1,2,\cdots,N\},\quad E_b=\{(i,j)| A_{I_iI_j} \text{ has nonzero entries} \}\]

\section{Hierarchical solver}
\label{sec:hsolver}

Although the Gaussian elimination works for any linear system, it destroys the sparsity when applied to sparse linear systems by introducing additional nonzero entries, or {\it fill-ins}, into the factorized matrices $L$ and $U$. Many attempts have been made to reduce these fill-ins via approximations; for example, incomplete LU factorization (ILU) strictly controls the number of fill-ins introduced.  
A key observation is that the fill-ins corresponding to well-separated interactions are low-rank \cite{Pouransari.H;Coulier.P;Darve.E2015a}.  For brevity, we simply use {\it h-solver} to denote this type of hierarchical solvers.  In the rest of this section, we introduce the algorithm in \cite{Pouransari.H;Coulier.P;Darve.E2015a} using block matrix notation, which is necessary for introducing the improved version of the algorithm in the next few sections.

\begin{figure}[htbp]
\begin{center}
\includegraphics[width=0.3\textwidth]{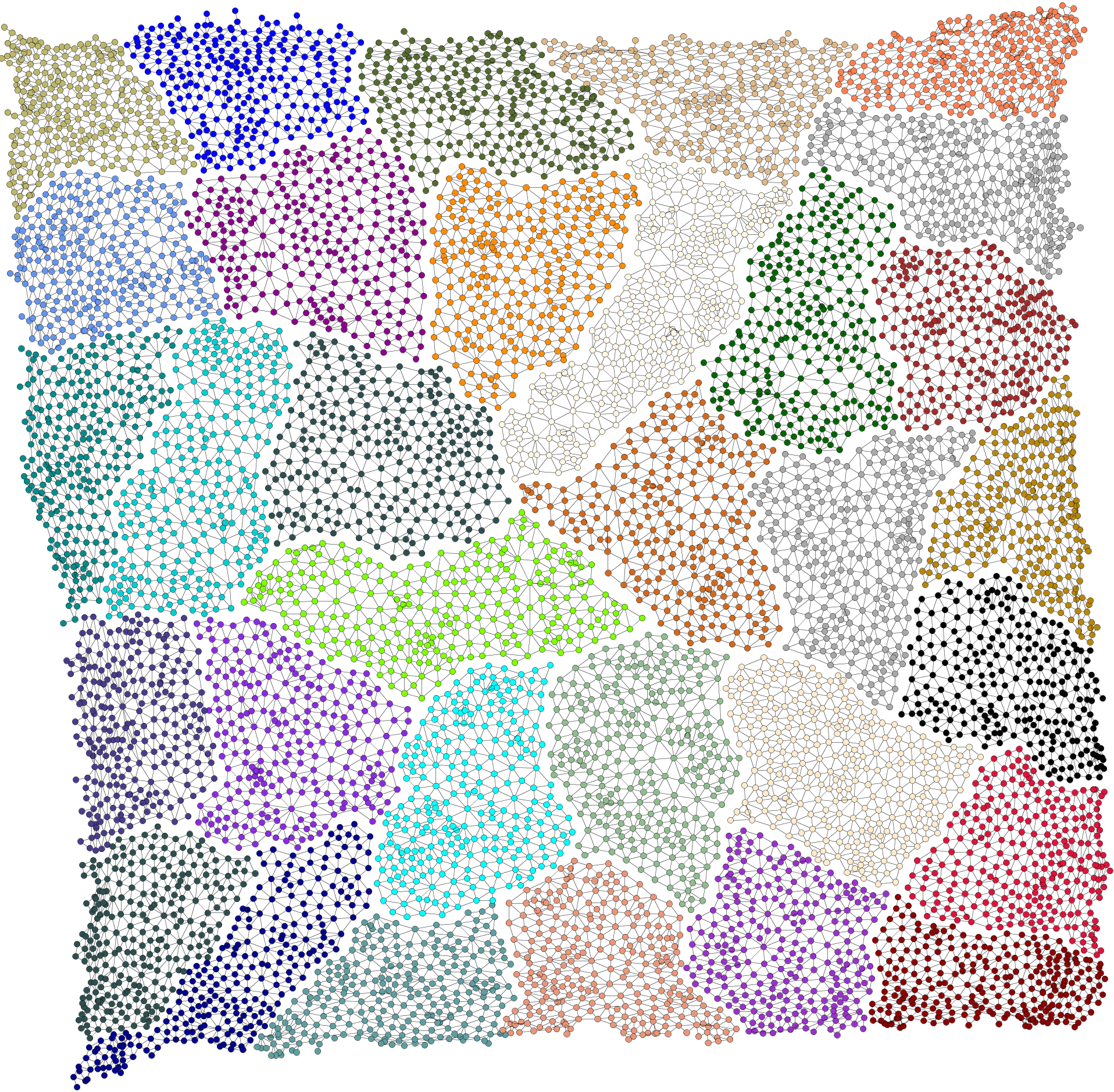}
\includegraphics[width=0.3\textwidth]{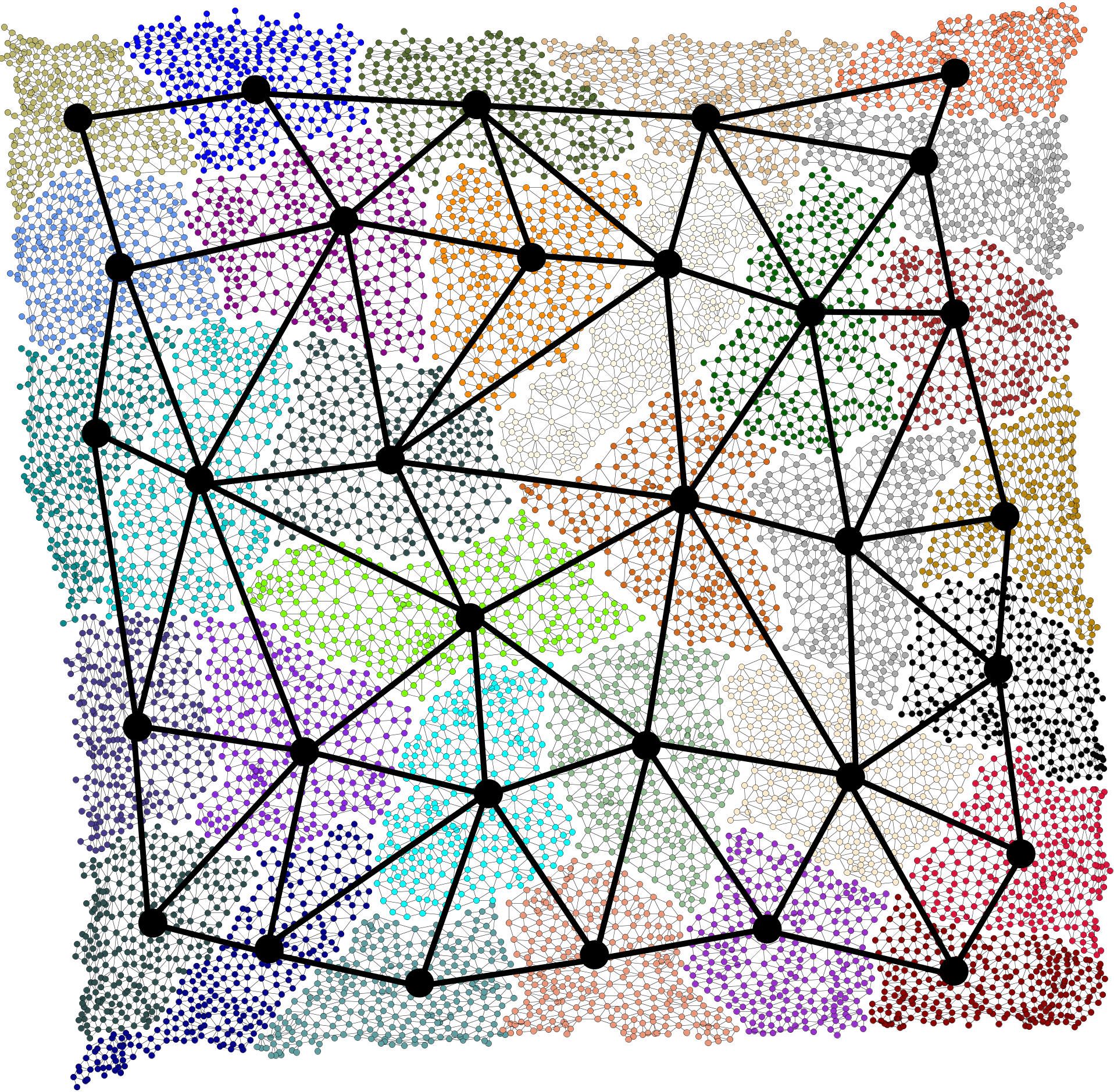}
\end{center}
\caption{Partitioning of the graph of a matrix resulting from a 2D discretization. Left: partitioning of the graph into local clusters denoted by different colors.  Right: The corresponding graph for the block version of the underlying matrix.
\label{fig:partition}}
\end{figure}

Given a sparse linear system, the first step of the proposed algorithm is to partition the unknowns into local clusters. The purpose of clustering is to accelerate the algorithm using BLAS3 operations \cite{Anderson.E;Bai.Z;Bischof.C;Blackford.S;Demmel.J;Dongarra.J;Du-Croz.J;Greenbaum.A;Hammerling.S;McKenney.A;others1999a}.  For example, the partitioning can be formed as is shown in Fig.~\ref{fig:partition}, where each color corresponds to one cluster.  This block partitioning corresponds to a new graph with fewer vertices.  An example of the new graph is given in Fig. \ref{fig:partition}.  As we will pairwisely agglomerate the clusters, we assume that there are $2N$ such clusters, denoted by $\{I_i^r\}_{i=1}^{2N}$.  Each cluster $I_i^r$ is call {\it red node} in \cite{Pouransari.H;Coulier.P;Darve.E2015a} and we continue to use this name to distinguish them from other types of clusters.  

Before we start the factorization, we first agglomerate the red nodes in a pairwise fashion. The motivation for this step is to facilitate the low-rank compression and will be explained in more details later in this paper. Let $\{I^s_i\}_{i=1}^{N}$ be the clusters formed by
$$I_i^s = I_{2i-1}^r\cup I_{2i}^r,\quad i=1,2,\ldots, N$$
The new clusters $\{I^s_i\}_{i=1}^{N}$ are called {\it super nodes} in \cite{Pouransari.H;Coulier.P;Darve.E2015a}.

\subsection{Adjacency and well-separated interactions}

Note that adjacency is based on the original graph before any elimination. During the elimination, new edges are introduced, but the adjacency between the super nodes remains unchanged throughout the elimination process.

The interactions between clusters are the off-diagonal blocks in the matrix. For example, the interaction between cluster $I_i$ and $I_j$ ($i\neq j$) is $A_{I_iI_j}$. This rule applies to both red nodes and super nodes and defines the adjacency. 

Although the adjacency stays the same, the interactions between clusters change during elimination.  Before we eliminate the first super node, the interactions are only between adjacent/neighboring clusters.  After we eliminate some clusters, interactions between non-adjacent clusters might be introduced into the matrix as fill-ins. 

The key idea of \cite{Pouransari.H;Coulier.P;Darve.E2015a} is closely related to the concept of {\bf well-separated interactions}, which is borrowed from the fast multipole method (FMM). In the FMM, two clusters are {\bf well-separated} if their distance is larger than a fraction of their diameters. A key finding in the FMM is that the interactions between well-separated clusters can be approximated accurately by low-rank matrices, assuming that their interactions are represented by non-oscillatory (``smooth'') kernels \cite{Greengard.L;Rokhlin.V1987a,Darve.E2000a,Ying.L;Biros.G;Zorin.D2004a}.  As the h-solver is purely algebraic, it does not have access to the geometric information as the FMM does.  Instead, well-separated interactions have to be defined purely based on the graph of the system matrix. One of the simplest approaches is to define that two clusters are well-separated if they are not neighbors.  Although different from the criterion in FMM, this type of well-separated interactions is expected to exhibit similar low-rank properties as those emerging from integral equations.  The low-rank properties have been numerically verified in \cite{Pouransari.H;Coulier.P;Darve.E2015a}.
\subsection{Factorization} 
The factorization phase of the h-solver is similar to that of a Gaussian (or in this instance Cholesky) elimination.  The major difference between an h-solver and a Gaussian elimination is that an h-solver first compresses well-separated interactions before it eliminates the current node.  

For each super node, an h-solver first compresses the well-separated interactions. Then it introduces additional variables and obtains an extended system matrix. After the extension, the current super node only has interactions with neighboring nodes; namely, there is no well-separated interaction between the current super node and the rest of the nodes. Then it performs the elimination. We will elaborate on the compression, extension, and elimination steps of the h-solver.  Assume that a given ordering of the super nodes is provided: $\{s^1,s^2,\ldots,s^N\}$.

Let's consider the case when we are eliminating a super node $s^i$.  We partition the unknowns of the linear system into three groups, corresponding to
\begin{itemize}
\item the current super node $s^i$, denoted by subscript ``s'',
\item the neighboring nodes of $s^i$, denoted by ``n'',
\item the well-separated nodes from $s^i$, denoted by ``w''.
\end{itemize}
Then the linear system (\ref{eq:Axb}) has the following block form
\begin{equation}
\label{eq:snw}
\begin{pmatrix}
A_{ss}&A_{sn}&A_{sw}\\
A_{ns}& A_{nn}&A_{nw}\\
A_{ws}&A_{wn}& A_{ww}\\
\end{pmatrix}
\begin{pmatrix}
x_s\\
x_n\\
x_w\\
\end{pmatrix}
=
\begin{pmatrix}
b_s\\
b_n\\
b_w
\end{pmatrix}
\end{equation}
If $s^i$ is not the first super node in the list, we consider that the super nodes that have already been eliminated are not included in (\ref{eq:snw}).  Moreover, new nodes will be introduced into the system in the extension step.  Therefore, among the neighboring and well-separated clusters of $s^i$, there are newly introduced parent-level red nodes.
We do not elaborate on these new red nodes here, but just assume the adjacency between $s^i$ and these new nodes are well defined.  We will define this adjacency in the extension step.

Note that when dealing with the first super node, there is no nonzero well-separated interactions; namely, $A_{sw}=0$.  $A_{sw}$ becomes nonzero during the algorithm due to fill-ins.

\subsubsection*{Compression}

The well-separated interactions $A_{sw}$ and $A_{ws}$ can be approximated well by low-rank matrices.  We perform a truncated Singular Value Decomposition (SVD)
\begin{equation}
\label{eq:local_compression}
A_{sw}=U\Sigma V^T+E_{sw}
\end{equation}
Here $U$ and $V$ are both orthonormal. $\Sigma$ is a diagonal matrix consisting of the most significant singular values of $A_{sw}$. The number of columns of $U$ is much smaller than the number of rows or columns of $A_{sw}$. Moreover, $E_{sw}$ is small and bounded by the given compression threshold.  For convenience, we define $R^T = \Sigma V^T$. Thus,
\begin{equation}
\label{eq:compression_sw}
A_{sw}=UR^T+E_{sw}
\end{equation}
We now solve the following approximate linear system
\begin{equation}
\label{eq:compressed}
\begin{pmatrix}
A_{ss}&A_{sn}&UR^T\\
A_{ns}& A_{nn}&A_{nw}\\
RU^T&A_{wn}& A_{ww}\\
\end{pmatrix}
\begin{pmatrix}
x_s\\
x_n\\
x_w\\
\end{pmatrix}
=
\begin{pmatrix}
b_s\\
b_n\\
b_w\\
\end{pmatrix}
\end{equation}

\subsubsection*{Extension}  In order to avoid more and more fill-ins, two additional unknowns, $y_b$ and $y_r$, are introduced.  Then we can instead solve
\begin{equation}
\label{eq:extended}
\begin{pmatrix}
A_{ss}&A_{sn}&& U&\\
A_{ns}& A_{nn}&A_{nw}&\\
&A_{wn}& A_{ww}&&R\\
U^T&&&&-I\\
&&R^T&-I&\\
\end{pmatrix}
\begin{pmatrix}
x_s\\
x_n\\
x_w\\
y_b\\
y_r\\
\end{pmatrix}
=
\begin{pmatrix}
b_s\\
b_n\\
b_w\\
0\\
0\\
\end{pmatrix}
\end{equation}
as $x_s, x_n, x_w$ in \eqref{eq:extended} solves \eqref{eq:compressed}.  In \cite{Pouransari.H;Coulier.P;Darve.E2015a}, the unknowns $y_b$ and $y_r$ are considered as a black node and a red node.  We will explain the underlying meanings of these new nodes at the end of this section.  Note that the newly introduced red node $y_r$ will not be merged or eliminated until all of the super nodes $\{s^1, s^2,\ldots, s^N\}$ are eliminated.

The red node $r$ needs well-defined neighbor list in order to continue the algorithm.  In fact, $r$ inherits the neighbors of $s$, which may contain both super nodes and red nodes. 

\subsubsection*{Elimination} 
In the elimination step, we eliminate $x_s$ and $y_b$. First, the linear system is rearranged as follows
\begin{equation}
\label{eq:extended_rearranged}
\begin{pmatrix}
{ A_{ss}}&{U}&{A_{sn}}&&\\
{U^T}&&&&-I\\
{A_{ns}}&& A_{nn}&A_{nw}&\\
&&A_{wn}& A_{ww}&R\\
&-I&&R^T&\\
\end{pmatrix}
\begin{pmatrix}
x_s\\
y_b\\
x_n\\
x_w\\
y_r\\
\end{pmatrix}
=
\begin{pmatrix}
b_s\\
0\\
b_n\\
b_w\\
0\\
\end{pmatrix}
\end{equation}

Denote (\ref{eq:extended_rearranged}) by 
\[K\tilde x = \tilde b\]

We eliminate the current super node $x_s$ and the black node $y_b$: $K=LK_2L^T$ with
\begin{align*}
L=
\begin{pmatrix}
I&&&&\\
L_{bs}& I&&\\
L_{ns}&L_{nb}&I&&\\
&&&I&\\
&S^{-1}&&&I\\
\end{pmatrix}&\\
K_2 =
\begin{pmatrix}
{A_{ss}}&&&&\\
&{-S}&&&\\
&& \tilde{\tilde {A}}_{nn}&A_{nw}&L_{nb}\\
&&A_{wn}& A_{ww}&R\\
&&L_{bn}&R^T&S^{-1}\\
\end{pmatrix}
=&
\begin{pmatrix}
A_{ss}&&\\
&-S&\\
&&A^{(i+1)}
\end{pmatrix}
\end{align*}
Here 
\begin{align*}
L_{bs} & = U^TA_{ss}^{-1} & L_{ns} & = A_{ns}A_{ss}^{-1}\\
L_{nb} & = A_{ns}A_{ss}^{-1}US^{-1} & S & = U^TA_{ss}^{-1}U\\
\tilde{\tilde {A}}_{nn} & = A_{nn}-A_{ns}A_{ss}^{-1}A_{sn} + A_{nb}S^{-1}A_{bn}
\end{align*}
In order to solve 
\begin{equation}
\label{eq:eliminated}
K_2L^T\tilde x=L^{-1}\tilde b
\end{equation} 
we need to solve a smaller linear system $A^{(i+1)}x=y$. By picking the next super node $s^{i+1}$, this linear system is again in the form of (\ref{eq:snw}).  We can repeat the compression, extension, and elimination procedure.

Note that each time these three steps are repeated, a super node $s^i$ is eliminated and a red node, denoted by $r^i$, is introduced to the system.  The black node is eliminated right after being introduced to the system.  By induction, we can conclude that at the beginning of each round of compression-extension-elimination, the nodes in the system are
$$\{s^i, s^{i+1},\ldots, s^N, r^1,\ldots, r^{i-1}\}$$
By the end of each round, the nodes in the system are 
$$\{s^{i+1},\ldots,s^N,r^1,\ldots,r^{i}\}$$  
This process is repeated until all the super nodes on the current level have been eliminated, by which time the nodes in the system are $\{r^1,r^2,\ldots,r^N\}$.  We denote the resulting matrix by $A^{\text{parent}}$.  Then we can again merge in a pairwise fashion the red nodes in order to form super nodes, and repeat the compression-extension-elimination for all the new super nodes until there is only one red node in the system.  Then we factorize the matrix exactly.

Note that as we repeatedly agglomerate pairs of red nodes into super nodes, we assume that the block partitioning of the original linear system (\ref{eq:Axb}) forms $N=2^{L}$ red nodes, where $L$ is an integer.  The assumption can be easily satisfied, for example, by recursively bisecting of the graph of (\ref{eq:Axb}).

In \cite{Pouransari.H;Coulier.P;Darve.E2015a}, this factorization algorithm is explained using hierarchical trees. We briefly explain the connection here and refer interested readers to \cite{Pouransari.H;Coulier.P;Darve.E2015a} for more details. 

\begin{figure}[htbp]
\begin{center}
\includegraphics[width=0.4\textheight]{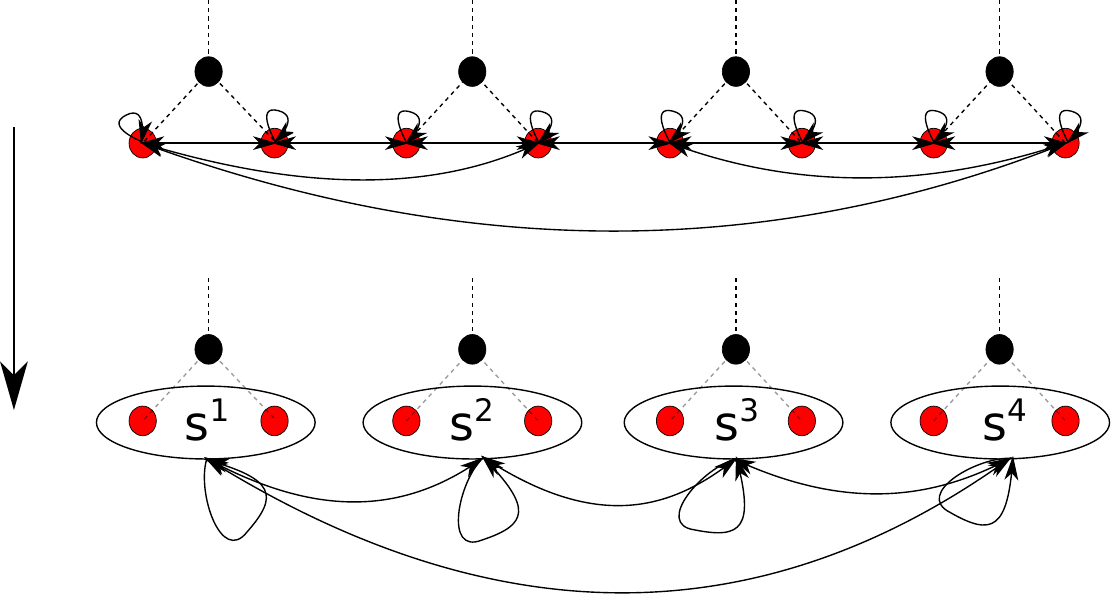}
\end{center}
\caption{Merge sibling red nodes \label{fig:merge}}
\end{figure}

\begin{figure}[htbp]
\begin{center}
\includegraphics[width=0.6\textwidth]{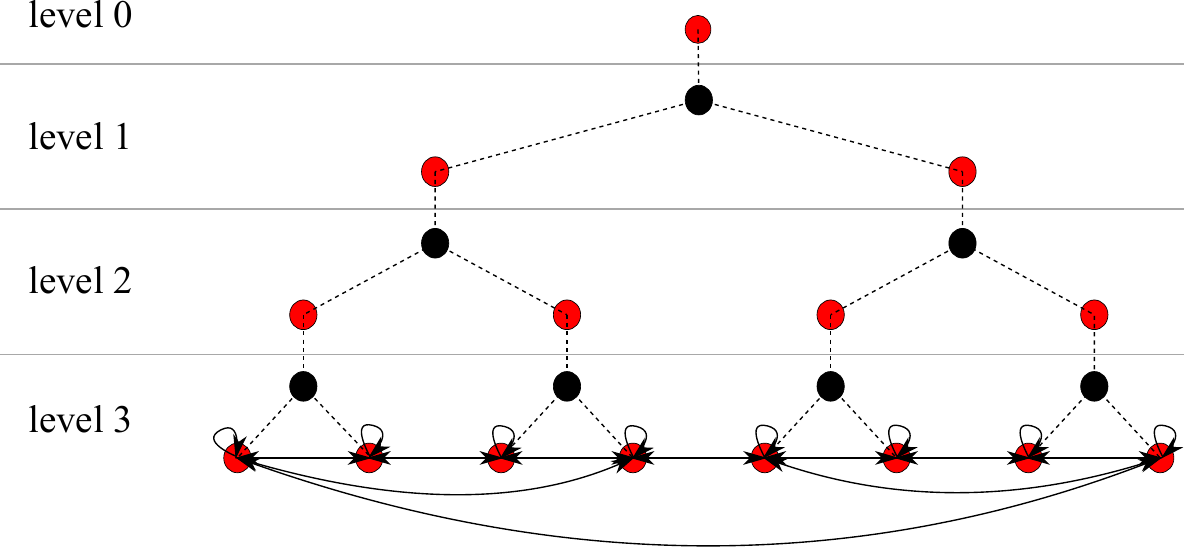}
\end{center}
\caption{An example of hierarchical tree
\label{fig:HTree}}
\end{figure}

The red nodes are pairwisely merged to form super nodes (Fig. \ref{fig:merge}), and each super node generates a black node and a red node. This procedure reduces the number of red nodes by half.  We put each black node on top of the corresponding super node, and then put the new red node on the corresponding black node.  After eliminating all the super nodes and the black nodes, we only have the newly introduced red nodes, which are only half as many as the red nodes we started with.  As this procedure will be repeated until there is only one red node, all of the nodes actually form a hierarchical tree. (See Fig. \ref{fig:HTree}.)

The tree notation of \cite{Pouransari.H;Coulier.P;Darve.E2015a} provides a vivid interpretation of the purely algebraic algorithm of the h-solver.  We refer to \cite{Pouransari.H;Coulier.P;Darve.E2015a} for more details on the tree formulations.
\subsubsection*{Computational complexity of the factorization}
It is shown in \cite{Pouransari.H;Coulier.P;Darve.E2015a} that the factorization algorithm has complexity $\mathcal{O}(nd_{\lm}^2)$, as long as there is a constant $\alpha$ such that $d_{i}/d_{i+1}<\alpha <2^{1/3}$, $\forall 0\leq i\leq \lm-1$. Here $n$ is the number of unknowns and $d_i$ is the maximum number of unknowns in a super node (or red node) on the $i$-th level.  The $\lm$-th level is the leaf level. The ratio $d_{i}/d_{i+1}$ depends on the low-rank property of the well-separated interactions.  For example, given the red node size $d_{i+1}$ on level $i+1$, the super node has size $2d_{i+1}$.  Suppose that the truncated SVD is able to provide the required accuracy with rank less that $d_{i+1}$, then the red node on level $i$ has size $d_i < d_{i+1}$. Then, $d_i/d_{i+1}$ satisfies the criterion.  We refer to \cite{Pouransari.H;Coulier.P;Darve.E2015a} for the details of complexity analysis.

\subsection{Solve}
After factorizing the matrix, we approximately solve $Ax=b$ with a given right hand side $b$.  This procedure can be denoted as solving the linear system $A_{\cH}x_{\cH}=b$.

Based on the factorization, the solution process is also recursive.  Starting with a right hand side $b$ and let $A^{(0)}=A$. If $i<N$, we eliminate a super node from $A^{(i)}$. For this, we split the right hand side $b$ as
\[
b=
\begin{pmatrix}
b_s\\
b_n\\
b_w\\
\end{pmatrix}
\]
We have to solve the extended system \eqref{eq:eliminated} recursively in three steps.
\begin{enumerate}
\item Forward substitution: $\phi=L^{-1}\tilde b$. Component-wise we have
\[
\phi=
\begin{pmatrix}
\phi_s\\
\phi_b\\
\phi_n\\
\phi_w\\
\phi_r\\
\end{pmatrix}
=
\begin{pmatrix}
b_s\\
-U^TA_{ss}^{-1}b_s\\
b_n-A_{ns}A_{ss}^{-1}b_s+A_{ns}A_{ss}^{-1}US^{-1}U^TA_{ss}^{-1}b_s\\
b_w\\
S^{-1}U^TA_{ss}^{-1}b_s\\
\end{pmatrix}
\]
\item Recursive solve: $\psi = K_2^{-1}\phi$. This step involves inverting $A^{(i+1)}$.  Component-wise we have:
\[
\psi=
\begin{pmatrix}
\psi_s\\
\psi_b\\
\psi_n\\
\psi_w\\
\psi_r\\
\end{pmatrix}
=
\begin{pmatrix}
A_{ss}^{-1}b_s\\
S^{-1}U^TA_{ss}^{-1}b_s\\
\psi_n\\
\psi_w\\
\psi_r\\
\end{pmatrix}
\quad\text{ and }\quad
A^{(i+1)}
\begin{pmatrix}
\psi_n\\
\psi_w\\
\psi_r\\
\end{pmatrix}
=
\begin{pmatrix}
\phi_n\\
\phi_w\\
\phi_r\\
\end{pmatrix}
\]
As for solving the second equation above, we go to the next super node and apply the algorithm in case $i+1< N$.  If $i+1=N$, we go to the next group of super nodes, which are on the next upper level in the hierarchical tree. 
\item Backward substitution: $\tilde x = L^{-T}\psi$. Component-wise we have
\[
\tilde x=
\begin{pmatrix}
x_s\\
y_b\\
x_n\\
x_w\\
y_r\\
\end{pmatrix}
=
\begin{pmatrix}
A_{ss}^{-1}(b_s-US^{-1}(U^TA_{ss}^{-1}(b_s-A_{sn}\psi_n)-\psi_r)-A_{sn}\psi_n)\\
S^{-1}(U^TA_{ss}^{-1}(b_s-A_{sn}\psi_n)-\psi_r)\\
\psi_n\\
\psi_w\\
\psi_r\\
\end{pmatrix}
\]
\end{enumerate}
Then the solution is
\[x=
\begin{pmatrix}
x_s\\
x_n\\
x_w\\
\end{pmatrix}
\]
This h-solver was named {\bf LoRaSp} ({ {\bf Lo}w {\bf Ra}nk {\bf Sp}arse solver}) in \cite{Pouransari.H;Coulier.P;Darve.E2015a}.

\section{Error Estimate}
\label{sec:error}
In this section, we consider how the local truncated SVDs contribute to the error of $A_{\cH}$. Note that in this paper we discuss the truncated SVDs, but other methods such as rank-revealing QR, rank-revealing LU, ACA could also be used.

 We will only consider a two-level method; namely, on the parent level the linear system is solved exactly.  There are $N$ steps of compression and elimination, where $N$ is the number of super nodes on the leaf level. 
 
We first recognize that the extension step could be expressed as a matrix product.
\begin{equation}
\label{eq:extension_large}
\begin{split}
&
\begin{pmatrix}
A_{ss}&A_{sn}&UR^T& &\\
A_{ns}& A_{nn}&A_{nw}&\\
RU^T&A_{wn}& A_{ww}&&\\
&&&&-I\\
&&&-I&\\
\end{pmatrix}\\
&=
\begin{pmatrix}
I&&&&U\\
&I&&&\\
&&I&R&\\
&&&I&\\
&&&&I\\
\end{pmatrix}
\begin{pmatrix}
A_{ss}&A_{sn}&& U&\\
A_{ns}& A_{nn}&A_{nw}&\\
&A_{wn}& A_{ww}&&R\\
U^T&&&&-I\\
&&R^T&-I&\\
\end{pmatrix}
\begin{pmatrix}
I&&&&\\
&I&&&\\
&&I&&\\
&&R^T&I&\\
U^T&&&&I\\
\end{pmatrix}
\end{split}
\end{equation}

In order to avoid constantly changing the size of the linear system due to the extension steps, we consider the extended system which contains all of the super nodes, the black nodes and the parent red nodes.  The following equation is equivalent to (\ref{eq:Axb}) in the sense that $x_s$ solves (\ref{eq:Axb}).  Note that the solid lines in this equation separate the three different types of nodes.
\begin{equation}
\label{eq:K0}
\left(
\begin{array}{c|c|c}
A&&\\\hline
&&-I\\\hline
&-I&\\
\end{array}
\right)
\begin{pmatrix}
x_s\\
x_b\\
x_r\\
\end{pmatrix}
=
\begin{pmatrix}
b\\
0\\
0\\
\end{pmatrix}
\end{equation}
Denote the matrix in (\ref{eq:K0}) by $\mathcal{K}_1$.  Then for each super node, the compression, extension and elimination steps can be reformulated in the extended system.   Consider the $i$-th super node. Right before the compression, the matrix is denoted by $\mathcal{K}_i$.  

In the compression step, we decompose the well-separated interactions into the low-rank part and error part, as is shown in (\ref{eq:compression_sw}).  We embed the error part $E_{sw}$ and $E_{ws}$ in the extended system with all types of nodes and denote it by $\cE_i$.  Let $\cK_i^-$ denote the matrix $\cK_i$ with the well-separated interactions $A_{sw}$ and $A_{ws}$ replaced by low-rank matrices $UR^T$ and $RU^T$, respectively.  Thus, $\cK_i=\cE_i+\cK_i^-$.

For the extension step, we first introduce the matrix 
\[
\cU_i=
\left(
\begin{array}{c|c|c}
I&R_i&U_i\\\hline
&I&\\\hline
&&I\\
\end{array}
\right)
\]
where the only nonzero parts of $R_i$ and $U_i$ are $R$ and $U$, respectively.  For example, $R$ is located on the rows corresponding to the nodes well-separated from the $i$-th super node, and on the columns corresponding to the $i$-th black node.  Similarly, we can locate $U$ in $U_i$.   Then the extension step is denoted by $\cK_i^-=\cU_i\cK_i^E\cU_i^T.$  Here $\cK_i^E$ corresponds to the extended matrix in the second line of (\ref{eq:extension_large}).

For the elimination step, we first embed $L$ in the extended system, with $1$'s on the rest of the diagonal entries, and denote it by $\cL_i$.  Then the elimination is denoted by $\cK_i^E=\cL_i\cK_{i+1}\cL_i^T$.  

Therefore, the compression, extension, and elimination steps for the $i$-th super node is denoted by 
\[
\cK_i=\cE_i+\cU_i\cL_i\cK_{i+1}\cL_i^T\cU_i^T
\]
If we recursively apply this equation, the factors $\cU_i\cL_i$ and $\cL_i^T\cU_i^T$ are applied on the error $\cE_j$ for any $i<j\leq N$.  Note that $\cL_i\cE_{i+1}\cL_i^T=\cE_{i+1}$, because $\cE_{i+1}$ is zero in the rows and columns corresponding to the first $i$ super nodes and black nodes.  Multiplying $\cU_i$ with $\cE_{i+1}$ may create additional nonzeros in the matrix.  The reason is that $\cE_{i+1}$ may include nonzero interaction between the $(i+1)$-th super node and the $i$-th parent red node, which we denote by $E_{ps}$.  In fact, the additional error term introduced by this multiplication is $UE_{ps}$.  As $U$ is the orthonormal low-rank basis, the additional error satisfies $\|UE_{ps}\|_F= \|E_{ps}\|_F$.

Moreover, note that $\cL_{i-1}\cU_{i}\cE_{i+1}\cU_i^T\cL_{i-1}^T=\cU_{i}\cE_{i+1}\cU_i^T$.   Then we can obtain the following equation
\[
\cK_1=\sum_{i=1}^{N} \tilde \cE_i +\cU_1\cL_1\cdots \cU_N\cL_N\cK_{N+1}\cL_N^T\cU_N^T\cdots \cL_1^T\cU_1^T
\]
with 
\[
\tilde \cE_i = 
\left\{
\begin{aligned}
&\cU_1\cdots\cU_{i-1}\cE_i\cU_{i-1}^T\cdots\cU_1^T,&i\geq 2\\
&\cE_1,& i=1\\
\end{aligned}
\right.
\]
As each $\cU_i$ introduces an additional error for the $i$-th super node, the total error can be estimated
\[
\|\tilde\cE_i\|_F^2\leq 2 \|\cE_i\|_F^2
\]
Let $\cK=\cK_1$, $\cE=\sum_i\tilde \cE_i$, and $\cK_{\cH}=\cU_1\cL_1\cdots \cU_N\cL_N\cK_{N+1}\cL_N^T\cU_N^T\cdots \cL_1^T\cU_1^T$.   $\cK_{\cH}$ is extended hierarchical factorization, while $\cK$ is basically the original linear system, with auxiliary variables introduced.  Assume that we apply an absolute error tolerance $\|\cE_i\|_F\leq\epsilon$. Then we get the following estimate
\[\|\cK-\cK_{\cH}\|_F\leq (2N)^{1/2}\epsilon\]

Although the hierarchical solver is realized by solving the extended linear system
\begin{equation}
\label{eq:KH}
\cK_{\cH}
\begin{pmatrix}
x_s\\
x_b\\
x_r\\
\end{pmatrix}
=
\begin{pmatrix}
b\\
0\\
0\\
\end{pmatrix}
\end{equation}
it could also be formulated as a linear system with the same number of unknowns as (\ref{eq:approx}).   After all, the hierarchical solver is a linear mapping from a right hand side $b$ to a solution $x_{\cH}$.

Consider a general extended linear system:
\begin{equation}
\label{eq:extended_ABC}
\begin{pmatrix}
A_-&B^T\\
B&C
\end{pmatrix}
\begin{pmatrix}
x\\
y
\end{pmatrix}
=
\begin{pmatrix}
b\\
0
\end{pmatrix}
\end{equation}

If the $x$ in the solution of the extended system (\ref{eq:extended_ABC}) solves the original system (\ref{eq:Axb}) for any $b$, we say that (\ref{eq:extended_ABC}) is an {\bf equivalent extension} of $(\ref{eq:Axb})$.  The following lemma provides a necessary and sufficient condition for an equivalent extension.
\begin{lemma}
Assume that $A$ and $C$ are invertible.  Then (\ref{eq:extended_ABC}) is an equivalent extension of (\ref{eq:Axb}) if and only if 
\[A_--B^TC^{-1}B=A\]
holds.  
\end{lemma}

Introduce the block partitioning, where ``f'' denotes the super nodes and ``c'' denote the black nodes and the parent red nodes.  Thus, we have the following partitioning:
\[
\cK
=
\begin{pmatrix}
K_{ff}&K_{fc}\\
K_{cf}&K_{cc}\\
\end{pmatrix}
\quad
\text{ and }
\quad
\cE
=
\begin{pmatrix}
 E_{ff}& E_{fc}\\
E_{cf}&\\
\end{pmatrix}
\]

As (\ref{eq:K0}) is an equivalent extension of (\ref{eq:Axb}), and (\ref{eq:KH}) is an equivalent extension of (\ref{eq:approx}), the following equations hold based on Lemma \ref{lm:extension}
\begin{align*}
A&=K_{ff}+E_{ff}-(K_{fc}+E_{fc})K_{cc}^{-1}(K_{cf}+E_{cf})\\
A_{\cH}&=K_{ff}- K_{fc} K_{cc}^{-1} K_{cf}
\end{align*}

Therefore, we can quantify the error
\[A-A_{\cH}=E_{ff}-E_{fc} K_{cc}^{-1} K_{cf}-K_{fc}K_{cc}^{-1}E_{cf}-E_{fc} K_{cc}^{-1}E_{cf}\] 
and obtain the following error estimate.
\begin{theorem}
The following inequality holds
\[
\|A-A_{\cH}\|_F\leq \sqrt{1+2\|K_{cc}^{-1}K_{cf}\|_F^2} \;
\|\cE\|_F+\| K_{cc}^{-1}\|_F\|\cE\|_F^2
\]
Moreover, assume that we have $N$ leaf super nodes and that each compression based on (\ref{eq:local_compression}) satisfies that $\|E_{sw}\|\leq \epsilon$. Then
\[
\|A-A_{\cH}\|_F\leq (2N)^{1/2}\sqrt{1+2\|K_{cc}^{-1} K_{cf}\|_F^2} \;
\epsilon+2N\| K_{cc}^{-1}\|_F \; \epsilon^2
\]
\end{theorem}

\section{Convergence Analysis}\label{sec:convergence}

In this section, we analyze the convergence of the hierarchical solvers we have introduced.  Approximation properties in the form of \eqref{eq:approx} are assumed and the conclusion applies to a wide range of approximate factorization solvers.  In fact, \cite{Bebendorf.M;Bollhofer.M;Bratsch.M2013a} has a similar convergence analysis. They assume that the approximate operator preserves vectors that are close to the small eigenvectors. Then they prove that the condition number of the preconditioned system $\kappa(A_{\cH}^{-1}A)$ is independent of the problem size.  In this section, we present a similar but more general convergence analysis.

Assume that a given collection of subsets of $\mathbb{R}^n$: $\{S_i\}_{i=1}^N$ satisfies the following conditions.
\begin{itemize}
\item[({\bf A1})] $\{S_i\}$ form a direct sum of $\mathbb{R}^n$; namely, $\mathbb{R}^n=S_1\oplus S_2\oplus\cdots\oplus S_N$. In addition, $\{S_i\}$ are $A$-orthogonal; namely, give any $i\neq j$, $v_i^TAv_j=0$ holds for $\forall v_i\in S_i$, $v_j\in S_j$.
\item[({\bf A2})] $A_{\cH}$ approximates $A$ up to relative error $\epsilon_i$ on $S_i$, i.e. 
\[\|(A-A_\cH)v\|\leq \epsilon_i\|A\|\|v\|,\quad\forall v\in S_i\]
\end{itemize}
Note that ({\bf A1}) is possible when  when $S_i$'s are spanned by nonoverlapping subsets of eigenvectors of $A$.   As $A$ is SPD, we use $\lambda_\text{max}$, the largest eigenvalue of $A$, instead of $\|A\|$ in the upper bound of ({\bf A2}).

Note that we can get an upper bound of $\kappa(A_{\cH}^{-1}A)$ by finding constants $c_0$ and $c_1$ such that
\[
c_0(Av,v)\leq (A_{\cH}v,v)\leq c_1(Av,v), \quad
\forall v\in \mathbb{R}^n 
\]
We have the following estimates for $c_1$ and $c_2$:
\begin{gather}
c_1=\max_{v}\frac{(A_{\cH}v,v)}{(Av,v)}=1+\max_{v}\frac{((A_{\cH}-A)v,v)}{(Av,v)}\label{eq:upperbound}\\
c_0=\min_{v}\frac{(A_{\cH}v,v)}{(Av,v)}=1-\max_{v}\frac{((A-A_\cH)v,v)}{(Av,v)}\label{eq:lowerbound}
\end{gather}
Define $P_i=S_i\oplus\cdots \oplus S_N.$ Then for $v_i\in P_i$, we can decompose it as $v_i=v_{i+1}+\tilde v_i$, with $v_{i+1}\in P_{i+1}$ and $\tilde v_i\in S_i$.  
\begin{align*}
&\max_{v_i\in P_i}\frac{((A_\cH-A)v_i,v_i)}{(Av_i,v_i)}\\
&=\max_{\tilde v_i\in S_i, v_{i+1}\in P_{i+1}}\frac{((A_\cH-A)(v_{i+1}+\tilde v_i),v_{i+1}+\tilde v_i)}{(A(v_{i+1}+\tilde v_i),v_{i+1}+\tilde v_i)}\\
&\leq\max_{\tilde v_i\in S_i, v_{i+1}\in P_{i+1}}\frac{((A_\cH-A)v_{i+1},v_{i+1})}{(Av_{i+1},v_{i+1})+(A\tilde v_{i},\tilde v_{i})} + \max_{\tilde v_i\in S_i, v_{i+1}\in P_{i+1}}\frac{((A_\cH-A)(2v_{i+1}+\tilde v_i),\tilde v_i)}{(A(v_{i+1}+\tilde v_{i}),v_{i+1}+\tilde v_i)}\\
&\leq\max_{v_{i+1}\in P_{i+1}}\frac{((A_\cH-A)v_{i+1},v_{i+1})}{(Av_{i+1},v_{i+1})} + \max_{\tilde v_i\in S_i, v_{i+1}\in P_{i+1}}\frac{\|(A_\cH-A)\tilde v_i\|\|2v_{i+1}+\tilde v_i\|}{(A(v_{i+1}+\tilde v_{i}),v_{i+1}+\tilde v_i)}\\
&\leq\max_{v_{i+1}\in P_{i+1}}\frac{((A_\cH-A)v_{i+1},v_{i+1})}{(Av_{i+1},v_{i+1})} + \max_{\tilde v_i\in S_i, v_{i+1}\in P_{i+1}}\frac{\epsilon_i\lambda_\text{max}\|\tilde v_i\|\|2v_{i+1}+\tilde v_i\|}{(A(v_{i+1}+\tilde v_{i}),v_{i+1}+\tilde v_i)}
\end{align*}
Define \[\mu_i:=\min_{v_i\in P_i}\frac{(Av_i,v_i)}{(v_i,v_i)},\quad 1\leq i\leq N\]
$\mu_i$ measures the lowest frequency of vectors in $P_i$.  For example, if $P_i$ is spanned by the eigenvectors $e_i, e_{i+1},\cdots,e_N$, then $\mu_i=\lambda_i$.

In order to estimate 
\[ \max_{\tilde v_i\in S_i, v_{i+1}\in P_{i+1}}\frac{\|\tilde v_i\|\|2v_{i+1}+\tilde v_i\|}{(A(v_{i+1}+\tilde v_{i}),v_{i+1}+\tilde v_i)}
\]
we simply assume that neither $\|\tilde v_i\|$ nor $\|2v_{i+1}+\tilde v_i\|$ is zero.  Then we have the following estimate:
\begin{align*}
&\frac{(A(v_{i+1}+\tilde v_{i}),v_{i+1}+\tilde v_i)}{\|\tilde v_i\|\|2v_{i+1}+\tilde v_i\|} \\
&=\frac{(A(2v_{i+1}+\tilde v_i),2v_{i+1}+\tilde v_i)+3(A\tilde v_i,\tilde v_i)}{4\|\tilde v_i\|\|2v_{i+1}+\tilde v_i\|} \\
&=\frac{\|2v_{i+1}+\tilde v_i\|}{4\|\tilde v_i\|}\frac{(A(2v_{i+1}+\tilde v_i),2v_{i+1}+\tilde v_i)}{\|2v_{i+1}+\tilde v_i\|^2}+\frac{\|\tilde v_i\|}{4\|2v_{i+1}+\tilde v_i\|}\frac{3(A\tilde v_i,\tilde v_i)}{\|\tilde v_i\|^2}\\
&\geq \mu_i\left(\frac{\|2v_{i+1}+\tilde v_i\|}{4\|\tilde v_i\|}+\frac{3\|\tilde v_i\|}{4\|2v_{i+1}+\tilde v_i\|}\right)\\
&\geq \frac{\sqrt{3}}{2}\mu_i
\end{align*}
Therefore, the following inequality holds
\[
\max_{v_i\in P_i}\frac{((A_\cH-A)v_i,v_i)}{(Av_i,v_i)}\leq\max_{v_{i+1}\in P_{i+1}}\frac{((A_\cH-A)v_{i+1},v_{i+1})}{(Av_{i+1},v_{i+1})} + \frac{2\epsilon_i\lambda_\text{max}}{\sqrt{3}\mu_i}
\]
We can repeat this procedure for $1\leq i\leq N-1$ and get
\[
\max_{v_i\in P_i}\frac{((A_\cH-A)v_i,v_i)}{(Av_i,v_i)}\leq \frac{\epsilon_N\lambda_\text{max}}{\mu_N}+\sum_{i<N } \frac{2\epsilon_i\lambda_\text{max}}{\sqrt{3}\mu_i} 
\]
Note that for $A-A_\cH$, the same bound can be derived:
\[
\max_{v_i\in P_i}\frac{((A-A_\cH)v_i,v_i)}{(Av_i,v_i)}\leq \frac{\epsilon_N\lambda_\text{max}}{\mu_N}+\sum_{i<N } \frac{2\epsilon_i\lambda_\text{max}}{\sqrt{3}\mu_i}
\]
Based on (\ref{eq:upperbound}), we get 
\[
c_1=1+\frac{\epsilon_N\lambda_\text{max}}{\mu_N}+\frac{2}{\sqrt{3}}\sum_{i<N }\frac{\epsilon_i\lambda_\text{max}}{\mu_i}
\]
In addition, assume $\epsilon_i$'s are small enough such that
\begin{equation}
\label{eq:coercive}
\frac{\epsilon_N\lambda_\text{max}}{\mu_N}+\frac{2}{\sqrt{3}}\sum_{i<N } \frac{\epsilon_i\lambda_\text{max}}{\mu_i}<1
\end{equation}
Then we get
\[
c_0=1-\frac{\epsilon_N\lambda_\text{max}}{\mu_N}-\frac{2}{\sqrt{3}}\sum_{i<N } \frac{\epsilon_i\lambda_\text{max}}{\mu_i}
\]
We now have an upper bound for the condition number $\kappa(A_{\cH}^{-1}A)$, as detailed in the theorem below.

\begin{theorem}
\label{thm:all_eval}
Assume ({\bf A1}), ({\bf A2}) and (\ref{eq:coercive}). Then, the following inequality holds 
\[\kappa(A_{\cH}^{-1}A)\leq \frac{1+\epsilon_N\lambda_\text{max}/\mu_N+\frac{2}{\sqrt{3}}\sum_{i<N}\epsilon_i\lambda_\text{max}/\mu_i}{1-\epsilon_N\lambda_\text{max}/\mu_N-\frac{2}{\sqrt{3}}\sum_{i<N }\epsilon_i\lambda_\text{max}/\mu_i}\]
\end{theorem}
\begin{example}
Assume we have a uniform accuracy $\epsilon$ for all eigenvectors. Then it is the special case $\mathbb{R}^n=S_1$ and, therefore, $\mu_1 =\lambda_\text{min}$, the minimum eigenvalue of $A$. Then we have 
\[\kappa(A_\cH^{-1}A)\leq \frac{1+\epsilon\lambda_\text{max}/\lambda_\text{min}}{1-\epsilon\lambda_\text{max}/\lambda_\text{min}}=\frac{1+\epsilon\kappa(A)}{1-\epsilon\kappa(A)}\]
This upper bound indicates that, with a uniform accuracy $\epsilon$ for error components of all frequency, $\epsilon$ has to decrease as fast as $1/\kappa(A)$ in order to keep $\kappa(A_{\cH}^{-1}A)$ bounded independent of $\kappa(A)$.
\end{example}

\begin{example}
\label{ex:two_piece}
Let $v_1,v_2,\ldots,v_m$ be a few vectors and $S_1:=\text{span}\{v_i\}_{i=1}^m$. Let $S_2:=S_1^{\perp_{A}}$ be the A-orthogonal complement of $S_1$ in $\mathbb{R}^n$. Then,

\[\kappa(A_{\cH}^{-1}A)\leq \frac{1+\frac{2}{\sqrt{3}}\epsilon_1\kappa(A)+\epsilon_2\lambda_\text{max}/\mu_2}{1-\frac{2}{\sqrt{3}}\epsilon_1\kappa(A)-\epsilon_2\lambda_\text{max}/\mu_2}\]
with
\begin{equation}
\label{eq:raleigh}
\mu_2 = \min_{v\in S_2}\frac{(Av,v)}{(v,v)}
\end{equation}
In particular, we can pick $0<\eta<1$ in order to get a partitioning of small and large eigenvalues:
\[S_1=\text{span} \{e_i\}_{i: \lambda_i< \eta},\quad S_2=\text{span}\{e_i\}_{i:\lambda_i\geq\eta}\] 
Then $\mu_2\geq\eta$. Therefore,
\[\kappa(A_{\cH}^{-1}A)\leq \frac{1+\frac{2}{\sqrt{3}}\epsilon_1\kappa(A)+\epsilon_2\lambda_\text{max}/\eta}{1-\frac{2}{\sqrt{3}}\epsilon_1\kappa(A)-\epsilon_2\lambda_\text{max}/\eta}\]
This example implies a practical strategy of improving the original hierarchical solver.  We consider $S_1$ to be spanned by a few low-frequency eigenvectors and improve the accuracy on these eigenvectors.  For the rest of the eigenvectors in $S_2$ we maintain the accuracy of the original hierarchical solver algorithm. Note that in general with $\eta$ being constant, the dimension of $S_1$ increases with the dimension of $A$.  
\end{example}

\begin{example}

Consider eigenvectors $e_1,e_2,\ldots,e_n$ and the corresponding eigenvalues $\lambda_1$ $\leq \lambda_2$ $\leq \cdots\leq\lambda_n$. Assume that $S_i:=\text{span}\{e_{i}\}$, $1\leq i\leq n$. Then $P_i=\text{span}_{j\geq i}\{e_{j}\}$ and $\mu_i=\lambda_i$.  Therefore, the following inequality holds
 \[\kappa(A_{\cH}^{-1}A)\leq \frac{1+\frac{2}{\sqrt{3}}\sum_{i<n}\epsilon_i\lambda_\text{max}/\lambda_i+\epsilon_n\lambda_\text{max}/\lambda_n}{1-\frac{2}{\sqrt{3}}\sum_{i<n }\epsilon_i\lambda_\text{max}/\lambda_i -\epsilon_n\lambda_\text{max}/\lambda_n}\]
 This upper bound suggests that the accuracy on each eigenvector shall be proportional to the magnitude of the corresponding eigenvalue in order to be efficient. 
\end{example}

%
%

\section{Improved hierarchical solver with guaranteed convergence}
\label{sec:improved}
Based on the convergence analysis, we can see that, in order to guarantee the convergence, we need to improve the accuracy of the h-solver on low-frequency error components. There are two approaches for this purpose.  The first approach is to let the $A_{\cH}$ preserves the action of $A$ on a few low-frequency vectors, such that the h-solver is exact on these error components.  The second approach is to improve the accuracy of $A_{\cH}$ on a few low-frequency vectors.  Both of these approaches add additional bases in the truncated SVD.  In this section, we introduce the first approach. The second approach is covered in {Appendix A}.

Note that in \cite{Bebendorf.M;Bollhofer.M;Bratsch.M2013a}, modifications to the hierarchical Cholesky decomposition are presented for the solver to satisfy the side constraint.  Although we have the same goal in this paper, this is a novel derivation due to the differences in the hierarchical factorization algorithms.

Based on Example \ref{ex:two_piece}, being exact on a few eigenvectors, i.e., letting $\epsilon_1=0$, leads to the following estimate
\[\kappa(A_{\cH}^{-1}A)\leq \frac{1+2\epsilon_2\lambda_\text{max}/\eta}{1-2\epsilon_2\lambda_\text{max}/\eta}\]

In this case if we can in addition guarantee constant $\eta$ and the same accuracy $\epsilon_2$ for problems of different sizes, then $\kappa(A_{\cH}^{-1}A)$ will be uniformly bounded independent of problem sizes.  However, in general this would be difficult to achieve.  For example, for the 2D Poisson equation discretized by the finite difference method on uniform grid, the number of eigenvalues smaller than a constant $\eta>0$ grows linearly with the number of unknowns.  Nevertheless, it is still promising that preserving a few eigenvectors may improve the performance of the hierarchical solvers.

Next, we consider how a vector is preserved during compression, extension, and elimination.

\subsubsection*{Preserving a vector during compression} During the compression step, only two off-diagonal block matrices are modified. It is straightforward to see that the modified blocks need to preserve the action of the original blocks on the corresponding segments of the vector.  The following algorithm provides sufficient conditions for preserving vectors during the compression.

\begin{lemma}
\label{lm:compression}
Given
\[
	\mathcal A=
	\begin{pmatrix}
	S&B^T&C^T\\
	B&D&E^T\\
	C&E&F\\
	\end{pmatrix}
	\quad \text{and} \quad
	\phi=
	\begin{pmatrix}
	\phi_x\\
	\phi_y\\
	\phi_z\\
	\end{pmatrix}
\]
let $\tilde{\mathcal{A}}$ be an approximation of $\mathcal{A}$ by replacing $B$ and $B^T$ by $\tilde B$ and $\tilde B^T$, respectively.  Assume that $B^T\phi_y=\tilde B^T\phi_y$ and $B\phi_x=\tilde B\phi_x$ hold. Then $\mathcal{A}\phi=\tilde{\mathcal{A}}\phi$.
\end{lemma}

Based on Lemma \ref{lm:compression}, the compression (\ref{eq:compression_sw}) has to satisfy
\begin{equation}
\label{eq:compression_xy}
A_{sw}\phi_y = UR^T\phi_y\quad\text{and }\quad A_{ws}\phi_x = RU^T\phi_x
\end{equation}
In order to satisfy this property, the compression of $A_{sw}$ has to follow a more complicated approach, instead of a straightforward SVD.  In general, we are looking for a decomposition of the form
\[
A_{sw} = UR^T+\hat U\hat R^T= [U,\hat U] 
\begin{bmatrix}
R^T\\
\hat R^T\\
\end{bmatrix}
\]
where the concatenated matrix $[U,\hat U]$ is orthonormal. 
  Then a sufficient condition for (\ref{eq:compression_xy}) to hold is 
\begin{equation}
	\phi_x\in \text{span}\{U\} \quad \text{and }\quad A_{sw}\phi_y\in \text{span}\{U\}
	\label{eq:preserve}	
\end{equation}
Therefore, we first orthonormalize the following matrix using QR or SVD 
\[[\phi_x,A_{sw}\phi_y] = U_1R^T_1\]
where the columns of $U_1$ are orthonormal.   We keep the component of $A_{sw}$ that is in the range of $U_1$, i.e., $U_1K_1^T$, with $K_1^T := U_1^TA_{sw}.$ Then, perform a low-rank factorization of the following matrix
\[(I-U_1U_1^T)A_{sw} = U_2K_2^T+E\]
Here $U_2K_2^T$ has low rank but approximates $(I-U_1U_1^T)A_{sw}$ well in the sense that $\|E\|/\|A_{sw}\|$ is less than the given compression parameter $\epsilon$. There are various approaches to obtain  this approximation, such as the truncated SVD, the rank-revealing QR and LU, ACA, etc. Then we obtain the following approximate factorization:
\[A_{sw}\approx [U_1,U_2] 
\begin{bmatrix}
K_1^T\\
K_2^T\\
\end{bmatrix}
\]
We can check that
\[ \Big \|
A_{sw}- [U_1,U_2] 
\begin{bmatrix}
K_1^T\\
K_2^T\\
\end{bmatrix}
\Big \|
\leq \epsilon \|A_{sw}\|
\]

\subsubsection*{Preserving a vector in the extension step}
During the extension step, the size of the linear system changes.  We need to consider carefully which vectors to preserve when the number of unknowns changes.  

Given the current global matrix $\cA$ and the vector $\phi$ we want to preserve, consider the following extended linear system
\[
\begin{pmatrix}
\cA_{-}&\cK^T\\
\cK&\cM\\
\end{pmatrix}
\begin{pmatrix}
x\\
y\\
\end{pmatrix}
=
\begin{pmatrix}
b\\
0\\
\end{pmatrix}
\]
the solution of which is consistent with that of $\cA x=b$ under the assumption $\cA=\cA_{-}-\cK^T\cM^{-1}\cK$.  For brevity we denote the extended system by $\cA_Ex_E=b_E$.  When $\cA_E$ is further approximated by 
\[\tilde \cA_E:=
\begin{pmatrix}
\tilde \cA_{-}&\tilde \cK^T\\
\tilde \cK&\tilde \cM
\end{pmatrix}
\]
in future steps, we also get an approximation for $\cA$, which is $\tilde \cA:= \tilde \cA_{-}-\tilde \cK^T\tilde \cM^{-1}\tilde \cK.$ We need to determine the vector to be preserved by the extended system such that $\tilde \cA\phi=\cA\phi$.  To be specific, in order to have $\tilde\cA\phi = \cA\phi$, we need to find a longer vector $\phi_E$ based on $\phi$ such that $\tilde\cA_E\phi_E=\cA_E\phi_E$ is a sufficient condition for $\tilde \cA\phi=\cA\phi$ to hold.

The following lemma provides us a guideline.

\begin{lemma}
\label{lm:extension}
Given $\phi$, let $\xi =-\cM^{-1}\cK\phi$ and $\phi_E=\binom{\phi}{\xi}$. Then a sufficient condition for $\tilde \cA\phi = \cA\phi$ to hold is $\tilde \cA_E\phi_E=\cA_E\phi_E$.
\end{lemma}

%
Although these equations seem complicated, it turns out that in the context of our algorithm, many of these terms simplify. In fact, the only calculation that needs to be done is $U^T \phi_x$, which is the representation (in the basis $U$) of the vector $\phi$ we wish to preserve for the current cluster $x$. $U^T \phi_x$ is the vector to preserve at the parent level (the red nodes right above the leaf level in the tree).

Here is a brief derivation of this result based on analyzing each step of the algorithm. Based on Lemma \ref{lm:extension}, in order to preserve 
\[
	\phi = 
	\begin{pmatrix}
			\phi_s\\
			\phi_n\\
			\phi_w
	\end{pmatrix}
\]
in the extension step we need to first calculate
\[
	\xi =
	-
	\begin{pmatrix}
		0&-I\\
		-I&0
	\end{pmatrix}^{-1}
	\begin{pmatrix}
		U^T & 0 & 0\\
		0 & 0 & R^T
	\end{pmatrix}
	\begin{pmatrix}
		\phi_s\\
		\phi_n\\
		\phi_w
	\end{pmatrix}
	=
	\begin{pmatrix}
		R^T\phi_w\\
		U^T\phi_s
	\end{pmatrix}
\]
The term $R^T\phi_w$ is associated with the black node, while $U^T\phi_s$ is for the parent red node. For the extended system $\cA_E$, the extended vector to be preserved is then
\[
	\phi_E =
	\begin{pmatrix}
		\phi_s\\
		R^T\phi_w\\
		\phi_n\\
		\phi_w\\
		U^T\phi_s
	\end{pmatrix}
\]

Now let's eliminate the leaf node and its black node, using a block Cholesky factorization.

\subsubsection*{Preserving a vector during elimination}

The elimination process can be written as
\[
	K = L K_2 L^T
\]
where $K_2$ is the partially factored matrix. Given an approximation $\tilde K_2$ of $K_2$, the matrix $\tilde K$
\[
	\tilde K:=L\tilde K_2L^T
\]
is our approximation of $K$.

\begin{lemma}
\label{lm:elimination}
Given a vector $\phi$, a sufficient condition for $K \phi=\tilde K\phi$ to hold is 
\[
	K_2\psi = \tilde K_2\psi \quad \text{where} \quad \psi =L^{T}\phi
\]
\end{lemma} 

The vector $\psi$ is of the form
\[
	\psi = L^{T} \phi
	=
	\begin{pmatrix}
		\tilde x_s\\
		\tilde y_b\\
		x_n\\
		x_w\\
		y_r
	\end{pmatrix} \quad \text{ given that }
	\quad
	\phi 
	=
	\begin{pmatrix}
		x_s\\
		y_b\\
		x_n\\
		x_w\\
		y_r
	\end{pmatrix}
\]
Here $\,\tilde{}\,$ denotes updated entries resulting from the multiplication by $L^{T}$. However, at that point in the algorithm the only blocks that need compression are those related to $x_n$, $x_w$ and $y_r$. The unknowns for $s$ and $b$ have already been eliminated so that $\tilde x_s$ and $\tilde y_b$ are automatically preserved in the rest of  the factorization.

In summary, looking back at all the steps, we simply proceed as follows. When we need to compress well-separated blocks, we calculate the vectors to preserve: $\phi_x$ and $A_{sw}\phi_y$ (see~\eqref{eq:preserve}) and then form the appropriate $U$ and $R$ matrices. We calculate $U^T \phi_x$, which becomes the vector associated with the parent red node and needs to be preserved during subsequent block low-rank approximations.

We name this improved version of h-solver {\bf GC}, which is short for {\it Guaranteed Convergence}.

\subsubsection*{Computational complexity of the GC}

Although in the GC we incorporate additional bases in the compression step, the criterion for linear complexity is still $d_i/d_{i+1}\leq\alpha<2^{1/3}$, $\forall 0\leq i\leq \lm-1$.  Note that $d_i$ consists of both the additional vectors incorporated and the bases from the truncated SVD.

\section{Numerical tests}
\label{sec:numerics}
In this section, we present numerical tests to demonstrate the performance of the h-solver and validate the theories proposed in previous sections.  The implementation is based on the LoRaSp\footnote{{\bf Lo}w {\bf Ra}nk {\bf Sp}arse solver, bitbucket.org/hadip/lorasp.} package.  LoRaSp depends on Eigen \cite{Guennebaud.G;Jacob.B;others2010a} for basic linear algebra calculations and SCOTCH \cite{Pellegrini.F;Roman.J1996a} for graph partitioning.  In the numerical tests, we solve the d-dimensional Poisson equation (d=2, 3)
\[
\left\{
\begin{aligned}
-\nabla\cdot(\alpha\nabla u)&=f\quad\text{in } \Omega:=[0,1]^d\\
u&=0\quad\text{on } \partial\Omega\\
\end{aligned}
\right.
\] 
discretized by the second-order central finite difference method. We solve both 2D and 3D Poisson equations. As we use standard five-point (2D) or seven point (3D) stencil, we can make use of the geometry information and define two clusters as well-separated clusters if their distance is larger than the size of the clusters. 

In the following numerical tests, we first compare the performance of the improved h-solver (GC) with that of the original h-solver (LoRaSp) for a constant-coefficient Poisson equation, i.e., $\alpha(x)=1$ in $\Omega$.   Note that as we are using the standard five-point (2D) or seven-point (3D) stencil, the eigenvalues and eigenvectors of the system matrices are available.  Therefore, we are able to enrich the GC with the smallest eigenvector.  In practice, computing eigenvectors might be computationally expensive.  Therefore, we also enrich the hierarchical basis with the constant vector, which approximates the low-frequency vectors well for the Poisson equation.  

In the second part of the numerical tests, we also solve a variable-coefficient Poisson equation. As the eigenvectors do not have explicit form for these cases, we only compare the LoRaSp with the GC-constant.

\subsection{Constant-coefficient Poisson equation}

\subsubsection*{2D cases}
In the 2D test, we first compare h-solvers by using them in the stationary iteration.  Given an initial guess $x_0$, the iterative scheme $x_{k+1}=x_k+A_{\cH}^{-1}(b-Ax_k)$ produces a sequence $x_0, x_1,\ldots$  As we randomly set up the solution $x^*$ and then obtain the right hand side $b=Ax^*$, the error of the $k$-th iteration can be evaluated as $e_k=x_k-x^*$. The iteration converges when the relative error $\|e_k\|/\|e_0\|$ is less than $10^{-6}$.

In \autoref{tb:2Dstationary}, we tested a sequence of problems with increasing number of unknowns $n$. The tree depths are adjusted accordingly to let the size of each leaf red node be $8$.  \autoref{tb:2Dstationary} shows the number of the stationary iteration for the relative error to decrease to less than $10^{-6}$.  In addition to the LoRaSp, we also tested the new h-solver enriched with a constant vector (GC-constant) and with the smallest eigenvector (GC-eigenvector). 


\begin{table}[htbp]
\caption{Number of the stationary iteration of the h-solvers for the constant coefficient Poisson equation. Leaf red nodes have a size of 8. The stopping criterion is relative error less than $10^{-6}$. The compression parameter $\epsilon$ is $0.1$ for all cases.}
\begin{center}
\begin{tabular}{ccccc}
\toprule
$n$ & tree depth & LoRaSp& GC-constant & GC-eigenvector  \\\midrule
$(2^5)^2$	 & 7	&4	&3	&3\\
$(2^6)^2$	 &9	    &5	&4	&4\\
$(2^7)^2$	 &11	&12	&5	&4\\
$(2^8)^2$	 &13	&24	&4	&5\\
$(2^9)^2$	 &15	&79	&5	&5\\
$(2^{10})^2$ &17	&289	&5	&6\\
\bottomrule
\end{tabular}
\end{center}
\label{tb:2Dstationary}
\end{table}%

From \autoref{tb:2Dstationary}, the number of iterations of the LoRaSp increases as the number of unknowns grows. In comparison, the GC-constant and the GC-eigenvector have an almost constant number of iterations for all cases. 

In practice, preconditioned Krylov subspace iterations are much more robust than the stationary iteration.  In \autoref{tb:2Dgmres}, we show the number of GMRES iterations preconditioned by the h-solvers. For each of the h-solvers, compression parameters $\epsilon = 0.1, 0.2,$ and $0.3$ are tested.  


\begin{table}[htbp]
\begin{center}
\begin{tabular}{ccccccccccc}
\midrule
& tree depth & \multicolumn{3}{c}{LoRaSp} & \multicolumn{3}{c}{GC-constant} & \multicolumn{3}{c}{GC-eigenvector} \\\midrule
$n \backslash \epsilon$	&	& $0.1$  & $0.2$& $0.3$& $0.1$  & $0.2$& $0.3$& $0.1$  & $0.2$& $0.3$\\\midrule
$(2^5)^2$		& 7	&5	&6	&7	&5	&6	&7	&5	&6	&7\\
$(2^6)^2$		&9	&6	&8	&10	&6	&7	&8	&6	&7	&8\\
$(2^7)^2$		&11	&8	&11	&14	&7	&8	&10	&6	&8	&10\\
$(2^8)^2$		&13	&11	&17	&22	&7	&9	&11	&7	&9	&11\\
$(2^9)^2$		&15	&18	&30	&39	&7	&10	&14	&8	&10	&14\\
$(2^{10})^2$	&17	&31	&53	&71	&8	&11	&16	&8	&11	&16\\
\bottomrule
\end{tabular}
\end{center}
\caption{Number of GMRES preconditioned by h-solvers solving constant coefficient Poisson equation. The tolerance $\epsilon$ used in the h-solver is shown on the second row. Leaf red nodes have a size of 8. The stopping criterion is relative residual less than $10^{-10}$. \label{tb:2Dgmres}}
\end{table}

For all of the h-solvers, increasing $\epsilon$ results in a larger number of iterations. For each value of $\epsilon$,  the LoRaSp shows a significantly increasing number of iterations.  The number of iterations of the GC-constant and the GC-eigenvector, however, are almost constant and only increase mildly with the number of unknowns. We used only relatively large values for $\epsilon$ since this leads already to a small number of iterations. $\epsilon$ below 0.1 reduces the number of iterations further. However, it leads to a more expensive preconditioner so that overall there is no gain.

In both the stationary iteration and preconditioned GMRES, the performance of the GC-constant and the GC-eigenvector are almost identical. This indicates that the constant vector is a good alternative for the smallest eigenvector for the 2D Poisson equation.

\subsubsection*{3D cases}

For the 3D Poisson equation, the condition number of the system matrices grows less rapidly with respect to the number of unknowns than in the 2D cases, although the storage cost and CPU time per iteration become higher due to denser system matrices.

In \autoref{tb:3Dgmres}, we show the number of GMRES iterations when preconditioned by the LoRaSp and the GC with $\epsilon=0.2$ and $0.3$. Moreover, as in \cite{Bebendorf.M;Bollhofer.M;Bratsch.M2016a}, an adaptive choice of $\epsilon$ for different problem sizes is proposed. We also apply this adaptive strategy in order to explore its influence on the h-solvers. 

In \cite{Bebendorf.M;Bollhofer.M;Bratsch.M2016a}, four approaches are tested:
\begin{itemize}
\item $\mathcal{H}\text{Chol}$: constant compression parameter for all cases. It is comparable to the LoRaSp with constant $\epsilon$.
\item $\mathcal{H}\text{Chols}$: adaptive compression parameter $\epsilon_l=\epsilon hd_l$, where $d_l$ is the minimum diameter of the finite element basis functions on level $l$. It is comparable to the LoRaSp with adaptive $\epsilon$. As we use finite difference and regular partitioning, we have $d_l=2^{(\lm-l)/3}h$.  If we use this criterion, $\epsilon_l$ becomes too small and is clearly sub-optimal for our method.  Therefore, we adopt a different compression parameter $\epsilon_l=\epsilon h2^{(\lm-l)/3}$ (we drop a factor $h$) such that the computation for large problems is more affordable. With this criterion, the relative accuracy between leaves and root is the same as in \cite{Bebendorf.M;Bollhofer.M;Bratsch.M2016a} (large $\epsilon$ at the root), but the scaling with $h$ is $O(h)$ instead of $O(h^2)$.
\item M$\mathcal{H}$Chol: preserving the constant vector and constant compression parameter. This is comparable to the GC-constant. 
\item M$\mathcal{H}$Chols: preserving the constant vector and adaptive compression parameter $\epsilon_l = \epsilon hD_l^{-1}$ ($D_l$: maximum diameter at level $l$). This is comparable to the GC with adaptive $\epsilon_l=\epsilon 2^{(l-\lm)/3}$ (small $\epsilon$ at the root). The compression parameter for the leaf level remains constant when the number of unknowns changes.  
\end{itemize}

\begin{table}[htbp]
\begin{center}
\begin{tabular}{ccccccccc}
\toprule
& tree depth      &   \multicolumn{3}{c}{LoRaSp } & \multicolumn{4}{c}{GC-constant}\\\midrule
$n\backslash\epsilon$
	 &      &  $0.2$ & $0.3$ & $\epsilon_{\lm}=8/\sqrt[3]{n}$& $0.2$ & $0.3$ &$\epsilon_{\lm}=0.2$& $\epsilon_{\lm}=0.3$\\\midrule
$8^3$	&6	& 5	&	5	&	7 		&5	&5	&5	&5\\
$16^3$	&9	&6	&	6	&	8 		&5	&6	&5	&6\\
$32^3$	&12	&7	&	8	&	9 		&6	&7 	&6	&7\\
$64^3$	&15	&8	&	11	&	9 		&6	&8	&6	&8\\
\bottomrule
\end{tabular}
\end{center}
\caption{Number of GMRES iterations with different hierarchical preconditioners solving constant coefficient Poisson equation. Leaf red nodes have a size of 8. The stopping criterion is a relative residual less than $10^{-10}$. In the LoRaSp, $\epsilon_{\lm}$ means that we use a small $\epsilon_l$ at the leaf and a constant $\epsilon_l$ at the root (see $\mathcal{H}\text{Chols}$). In the GC, $\epsilon_l$ is constant at the leaf (as indicated) and goes down to $O(h)$ at the root (see M$\mathcal{H}$Chols).
\label{tb:3Dgmres}}
\end{table}%

From \autoref{tb:3Dgmres}, the LoRaSp has a number of iterations growing very slowly with the number of unknowns, even though large values of $\epsilon$ (0.2, 0.3) are used. This is quite different from \cite{Bebendorf.M;Bollhofer.M;Bratsch.M2016a} where they report a much faster increase in the number of iterations with problem size. 
The adaptive choice of $\epsilon$ leads to a roughly constant number of iterations. We note that the adaptive $\epsilon_{\lm}$ keeps decreasing as $n$ increases, so that the cost of the preconditioner goes up with $n$ at a faster rate than $n$. If we assume that the rank is $O(\log 1/\epsilon)$, then the CPU time grows roughly like $O(n \log^2 n)$. A similar increase in CPU time occurs in \cite{Bebendorf.M;Bollhofer.M;Bratsch.M2016a}.

The GC-constant improves over the LoRaSp. In this case, the adaptive strategy leads to larger ranks near the root. However, since the number of iterations is already nearly constant, this choice for $\epsilon_l$ does not seem to make a difference. The CPU time is larger with adaptive $\epsilon_l$ compared to non-adaptive. So the adaptive choice does not seem preferable for this benchmark.

Different results might be obtained on larger grids. However, limitations in memory do not allow us to run larger test cases. The largest case we ran is $n \approx $ 262k.
  

\subsection{Variable-coefficient Poisson equation}
Variable-coefficient Poisson equation is frequently encountered in applications of engineering problems, for example, the heat conduction problem with more than one medium. Another example is fluid flow in heterogeneous porous media.  We consider the following two variable-coefficient cases.

\begin{figure}[htbp]
\setlength{\unitlength}{0.12in} 
\centering
\begin{picture}(10,10) 
\put(0,0){\line(0,1){10}}
\put(10,0){\line(0,1){10}}
\put(0,10){\line(1,0){10}}
\put(0,0){\line(1,0){10}}
\put(2.5,2.5){\line(0,1){5}}
\put(7.5,2.5){\line(0,1){5}}
\put(2.5,7.5){\line(1,0){5}}
\put(2.5,2.5){\line(1,0){5}}
\put(5,5){$\Omega_1$}
\put(.5,2.5){$\Omega_2$}
\end{picture}
\caption{Subdomains for piecewise-constant coefficient Poisson equation}
 \label{fig:subdomains}
\end{figure}
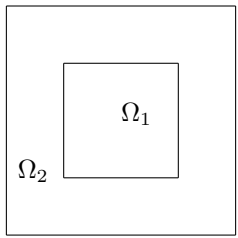
 
\begin{itemize}
\item {\bf Piecewise-constant coefficient}: Given $\Omega=[0,1]\times[0,1]$, define a subdomain $\Omega_1:=[1/4,3/4]\times[1/4,3/4]$. The rest of the domain is $\Omega_2:=\Omega\backslash \Omega_1$. (See Fig. \ref{fig:subdomains}.) Let 
\[\alpha(x)=\left\{
\begin{aligned}
&10^{-5},&\text{ if } x\in \Omega_1\\
&1,&\text{ if } x \in\Omega_2\\
\end{aligned}
\right.\]
\item {\bf Random coefficient}: we use a simplified model (not physical) in which we use uniformly distributed random variables in $[0,1]$ for $\alpha$ in the discretized system.
\end{itemize}

\begin{table}[htbp]
\begin{center}
\begin{tabular}{cccccccc}
\toprule
	      & tree depth & \multicolumn{3}{c}{LoRaSp}  & \multicolumn{3}{c}{GC-constant} \\\midrule
$n\backslash\epsilon$	&	& $0.1$  & $0.2$ & $0.3$ & $0.1$& $0.2$ & $0.3$\\\midrule
$(2^5)^2$		&7	&6	&6	&6	&7	&8	&7\\
$(2^6)^2$		&9	&6	&7	&8	&7	&8	&9\\
$(2^7)^2$		&11	&9	&10	&14	&8	&10	&12\\
$(2^8)^2$		&13	&11	&16	&21	&10	&10	&13\\
$(2^9)^2$		&15	&17	&29	&37	&9	&11	&15\\
$(2^{10})^2$	&17	&31	&54	&70	&10	&13	&18\\
\bottomrule
\end{tabular}
\end{center}
\caption{Number of GMRES iterations preconditioned by h-solvers for the piecewise-constant coefficient Poisson equation. Leaf red nodes have a size of 8. The stopping criterion is relative residual less than $10^{-10}$.
\label{tb:pwconstant}}
\end{table}%

\begin{table}[htbp]
\caption{Number of GMRES iterations preconditioned by h-solvers for the random coefficient Poisson equation. Leaf red nodes have a size of 8. The stopping criterion is relative residual less than $10^{-10}$.}
\begin{center}
\begin{tabular}{cccccccc}
\toprule
	      & tree depth & \multicolumn{3}{c}{LoRaSp}  & \multicolumn{3}{c}{GC-constant} \\\midrule
$n\backslash\epsilon$	&	& $0.1$  & $0.2$ & $0.3$ & $0.1$& $0.2$ & $0.3$\\\midrule
$(2^5)^2$		&7	&5	&6	&7	&5	&7	&7\\
$(2^6)^2$		&9	&7	&8	&9	&6	&7	&9\\
$(2^7)^2$		&11	&9	&11	&14	&7	&9	&11\\
$(2^8)^2$		&13	&12	&18	&22	&7	&10	&12\\
$(2^9)^2$		&15	&20	&31	&40	&8	&10	&15\\
$(2^{10})^2$	&17	&33	&56	&74	&8	&12	&18\\
\bottomrule
\end{tabular}
\end{center}
\label{tb:random}
\end{table}%

In \autoref{tb:pwconstant}, the number of GMRES iterations is shown for different hierarchical solvers for the piecewise-constant coefficient Poisson equation. In \autoref{tb:random}, the results for the random coefficient Poisson equation are shown. Both of these two tables show similar performance of the LoRaSp and the GC-constant compared with the constant coefficient case. These results show that the GC-constant is consistently improving over the LoRaSp. This indicates that the GC-constant is a promising algorithm for practical engineering problems. The GC-constant leads to a number of iterations that appears bounded (in our benchmarks), implying that the preconditioned system has a bounded condition number and that the eigenvalues have a favorable distribution for GMRES.

\section{Conclusion}

In this paper, we modified the algorithm of the hierarchical solver in \cite{Pouransari.H;Coulier.P;Darve.E2015a} and conducted a convergence analysis.  Based on the convergence analysis, we proposed an improved hierarchical solver that increases the accuracy on the low-frequency components of the error.  Given certain eigenvectors or an approximation thereof, the algorithm is able to exactly preserve these vectors. 

As a practical choice, we considered the use of piecewise constant vectors as local approximations of the smooth eigenvectors associated with small eigenvalues. This choice applies to elliptic PDEs discretizations. This is a computationally cheap option that provides a sufficiently good approximation of the eigenvectors. Under certain assumptions, this improved solver is an optimal preconditioner. In our benchmarks for example, we observed that GMRES converges with an almost constant number of iterations, regardless of the problem size.


\section*{Acknowledgments}

Funding from the ``Army High Performance Computing Research Center'' (AHPCRC), sponsored by the U.S.\ Army Research Laboratory under contract No.\ W911NF-07-2-0027, at Stanford, supported in part this research. The financial support is gratefully acknowledged.

\newpage
\begin{appendix}\label{ap:app}
\section{Detailed algorithms}
In this section we provide detailed algorithms for the setup and solve phase of LoRaSp.  The setup phase is summarized in a recursive fashion by Algorithm \ref{alg:setup}.  We use $A^{\text{leaf}}$ to denote the matrix in (\ref{eq:Axb}).
\begin{algorithm}
\caption{Setup phase of the h-solver}
\label{alg:setup}
\begin{algorithmic}
\State \Call{Setup}{$A^{\text{leaf}}$}
\Function{Setup}{$A$}
\For{each pair of sibling red nodes $r_1$ and $r_2$ in the red node list}
\State{Merge $r_1$ and $r_2$ into a super node $s$. }
\EndFor
\For{$s$ in the super node list}
\State{Compress the well-separated interaction $A_{sn}, A_{ns}$. See (\ref{eq:compressed}).}
\State{Introduce a black node $b$ and a parent red node $r$. Extend the matrix. See (\ref{eq:extended}).}
\State{Eliminate the super node $s$ and the black node $b$. See (\ref{eq:eliminated}).}
\EndFor
\If{the parent level is the root}
\State{Factorize $A^{\text{parent}}$ directly.}
\Else
\State{Setup($A^{\text{parent}}$).}
\EndIf
\EndFunction
\end{algorithmic}
\end{algorithm}

For the solve phase of LoRaSp, we first introduce some more notation to simplify the expressions:\begin{align*}
P_{ss}:&=S^{-1}U^TA_{ss}^{-1}\\
M_{ss}:&=A_{ss}^{-1}-P_{ss}^TSP_{ss}
\end{align*}
Note that $S  = U^TA_{ss}^{-1}U$.  Then some of the expressions can be simplified:
\begin{align*}
b_n-A_{ns}A_{ss}^{-1}b_s+A_{ns}A_{ss}^{-1}US^{-1}U^TA_{ss}^{-1}b_s&=b_n-A_{ns}M_{ss}b_s\\
A_{ss}^{-1}(b_s-US^{-1}(U^TA_{ss}^{-1}(b_s-A_{sn}\psi_n)-\psi_r)-A_{sn}\psi_n)&=M_{ss}(b_s-A_{sn}\psi_n)+P_{ss}^T\psi_r)
\end{align*}
Note that $P_{ss}$ maps the unknowns of the current super node to the unknowns of the corresponding parent red node. $M_{ss}$ is closely related to $A_{ss}^{-1}$.  In fact, when $U$ vanishes, $M_{ss}$ is exactly $A_{ss}^{-1}$.  On the other hand, when $U$ is invertible, $M_{ss}=0$. $M_{ss}$ can be factored as:
\[
	M_{ss} = 
	\Big( I - A_{ss}^{-1} U S^{-1}
	U^T \Big ) A_{ss}^{-1}
	= P_\text{perp} A_{ss}^{-1}
\]
where $P_\text{perp}$ is an $A_{ss}$ orthogonal projection onto $\text{span}\{U\}^{\perp}$. 
Here $\text{span}\{U\}$ denote the space spanned by the columns of $U$.

The solve algorithm involves recursion from one super node to the next super node.  However, it is possible to rearrange the order of computation to first finish the computation for all of the super nodes on one level of the hierarchical tree and then recurse to the next upper level of the tree. The solve step is formulated as Algorithm \ref{alg:solve}.  It makes a recursive call when moving on to the next super node.  If $i=N$, then all of the current super nodes are eliminated and the system matrix $A^{(N)}$ corresponds to parent level red nodes. Then, call the solve algorithm recursively on the parent level.
\begin{algorithm}
\caption{Solve phase of h-solver}
\label{alg:solve}
\begin{algorithmic}
\State $x=$ \Call{Solve}{$A^\text{leaf},b^{\text{leaf}}$}
\Function{Solve}{$A,b$}
\State{For all parent level node $p$, set $r_p = 0$.}
\For{ $s$ in the super node list}
\State{
Forward substitution for the neighbors $n$ of $s$:
\begin{center}
$\phi_n= r_n-A_{ns}M_{ss}r_s$
\end{center}
}
\EndFor
\For{$r$ in the parent level red node list}
\State{Restriction:
\begin{center}
$\phi_r =\phi_r+ P_{ss}\phi_s$
\end{center}
}
\EndFor
\State{Parent level solve:
\begin{center}
$\psi^\text{parent}=$\Call{Solve}{$A^\text{parent},\phi^\text{parent}$}
\end{center}
}
\For{$r$ in the parent level red node list}
\State{Prolongation:
\begin{center}
$x_s = P_{ss}^T\psi_r$
\end{center}
}
\EndFor
\For{$s$ in the super node list of a reversed order}
\State{
Backward substitution for the neighbors $n$ of $s$:
\begin{center}$x_s = x_s+M_{ss}(b_s-A_{sn}x_n)$
\end{center}
}
\EndFor
\State \Return $x$
\EndFunction
\end{algorithmic}
\end{algorithm}
\section{Improved convergence on low-frequency eigenvectors}
In this section, we provide alternative approaches to only approximately preserve vectors for GC.  Based on \autoref{thm:all_eval}, we need the accuracy on each eigenvector to be ideally proportional to the magnitude of the corresponding eigenvalue.  In order to be computationally efficient, we only improve the accuracy for a few smallest eigenvectors.  In the rest of this appendix, we introduce three approaches for this purpose: {\bf first-order one-sided projection}, {\bf first-order symmetric projection}, and {\bf second-order symmetric projection}.  The compression is mostly related to the well-separated block $A_{sw}$ and $A_{ws}$.  For brevity, we denote $A_{ws}$ by $B$.

During one step of compression, the well-separated interactions $B$ and $B^T$ are approximated by low-rank matrices. In order to improve the convergence on a few low-frequency eigenvectors: $\{e_1,e_2,\ldots, e_M\}$, we need the approximations $\tilde B$ and $\tilde B^T$ to be more accurate on $\{[e_i]_s\}$ and $\{[e_i]_w\}$, respectively. Here $[\cdot]$ denotes a segment of a vector and subscript $s$ and $w$ denote the current super node and its well-separated nodes, respectively. Our objective is to find $\tilde B$ such that 
\begin{gather}
\label{eq:o2}
\begin{split}
\|(B-\tilde B)[e_i]_s\| \leq \epsilon_1 \lambda_i\|B\|\|[e_i]_s\|, \quad\forall e_i, \;\; {i\leq M} \\
\|(B^T-\tilde B^T)[e_i]_w\| \leq \epsilon_1\lambda_i \|B\|\|[e_i]_w\|, \quad\forall e_i, \;\; {i\leq M}
\end{split} \\
\label{eq:o3}
\|B-\tilde B\|\leq \epsilon_2\|B\|
\end{gather}
Note that $\epsilon_1\leq\epsilon_2$. The objectives mean that $\tilde B$ approximates $B$ with an accuracy of $\epsilon_2$. Moreover, the approximation has an improved accuracy of $\epsilon_1$ on $\text{span}\{e_i\}_{i\leq N}$, and the accuracy is proportional to $\lambda_i$ on the direction of $e_i$, $\forall 1\leq i\leq N$.  Let 
\[X_s=\left[\frac{[e_1]_s}{\lambda_1},\frac{[e_2]_s}{\lambda_2},\ldots,\frac{[e_M]_s}{\lambda_M}\right]\]
and
\[X_w=\left[\frac{[e_1]_w}{\lambda_1},\frac{[e_2]_w}{\lambda_2},\ldots,\frac{[e_M]_w}{\lambda_M}\right]\]

\subsection*{First-order one-sided projection}
We first perform a truncated SVD for $[X_s,B^TX_w]$ and obtain $U_1$, the purpose of which is to approximately preserve $[X_s,B^TX_w]$.  Then we keep the component of $[X_s,B^T]$ in the ``direction'' of $U_1$ and compress the rest of it via another truncated SVD. 
\begin{enumerate}
\item Perform a truncated SVD:
\[\left[X_s,B^TX_w\right] = U_1\Sigma_1 V_1^T+\hat U_1\hat\Sigma_1\hat V_1^T\]
where $\|\hat U_1\hat\Sigma_1\hat V_1^T\|\leq \epsilon_1\|[X_s,B^TX_w]\|$.
Then 
\[\|(I-U_1U_1^T)[X_s,B^TX_w]\|\leq \epsilon_1\|[X_s,B^TX_w]\|\]
Obtain an intermediate approximation
\[\hat B^T = U_1U_1^TB^T\]

\item Perform another truncated SVD:
\[(I-U_1U_1^T)B^T=U_2\Sigma_2V_2^T+E_2
\quad
\text{with}
\quad
\|E_2\|\leq \epsilon_2\|B^T\|\]
Obtain the final approximation 
\[\tilde B^T:=\hat B^T+U_2U_2^T(B^T-\hat B^T)\]
\end{enumerate}
Note that for any given vector $v_s$ and $v_w$, we alway have
\begin{align*}
\|(B-\tilde B)v_s\|&=\|B(I-U_1U_1^T)(I-U_2U_2^T)v_s\|\leq \epsilon_2\|B\|\|v_s\|\\
\|(B^T-\tilde B^T)v_w\|&=\|(I-U_2U_2^T)(I-U_1U_1^T)B^Tv_w\|\leq \epsilon_2\|B^T\|\|v_w\|
\end{align*}
Therefore,
\begin{align*}
\|(B-\tilde B)X_s\alpha_i\|&=\|B(I-U_1U_1^T)(I-U_2U_2^T)(I-U_1U_1^T)X_s\alpha_i\|\\
&\leq \|B(I-U_1U_1^T)(I-U_2U_2^T)\| \|(I-U_1U_1^T)X_s\|\|\alpha_i\|\\
&\leq \lambda_i\epsilon_1\epsilon_2\|B\|\|[X_s,B^TX_w]\|\|\alpha_i\|
\end{align*}

For $B^T-\tilde B^T$, we have the following estimates
\begin{align*}
\|(B^T-\tilde B^T)X_w\| & = \|(I-U_2U_2^T)(I-U_1U_1^T)B^TX_w\|\leq \epsilon_1\|[X_s,B^TX_w]\|\\
\|(B^T-\tilde B^T)X_w\| & = \|(I-U_2U_2^T)(I-U_1U_1^T)B^TX_w\|\leq \epsilon_2\|B^T\|\|X_w\|
\end{align*}
Therefore, 
\[\|(B^T-\tilde B^T)X_w\|\leq \min\{ \epsilon_1\|[X_s,B^TX_w]\|,\epsilon_2\|B^T\|\|X_w\|\}\]

\subsubsection*{First-order symmetric projection}
In this subsection we introduce a symmetric approach to approximate preserve vectors.  The approximation $\tilde B^T$ is obtained via the following steps.
\begin{enumerate}
\item We first perform two truncated SVDs:
\begin{align*}
 \left[X_s,B^TX_w\right] &=U_1^s\Sigma_1^s (V_1^s)^T+\hat U_1^s\hat\Sigma_1^s(\hat V_1^s)^T \\
[X_w, BX_s]&=U_1^w\Sigma_1^w (V_1^w)^T+\hat U_1^w\hat\Sigma_1^w(\hat V_1^w)^T
\end{align*}
which satisfy 
\[\|\hat U_1^s\hat\Sigma_1^s(\hat V_1^s)^T\|\leq \epsilon_1\|[X_s,B^TX_w]\| \qquad \|\hat U_1^w\hat\Sigma_1^w(\hat V_1^w)^T\|\leq \epsilon_1\|[X_w,BX_s]\|\]
Note that equivalently we have
\begin{align*}
\|(I-U_1^s(U_1^s)^T)[X_s,B^TX_w]\|&\leq \epsilon_1\|[X_s,B^TX_w]\|\\
 \|(I-U_1^w(U_1^w)^T)[X_w,BX_s]\|&\leq \epsilon_1\|[X_w,BX_s]\|
\end{align*}

Then an intermediate approximation is 
\[\hat B^T=U_1^s(U_1^s)^TB^TU_1^w(U_1^w)^T\]

\item Perform another SVD:
\[B^T-\hat B^T= U_2\Sigma_2 V_2^T+\hat U_2\hat\Sigma_2\hat V_2^T\]
where $\|\hat U_2\hat\Sigma_2\hat V_2^T\|\leq \epsilon_2 \|B^T\|,$ or, equivalently,  
\[\|(I-U_2U_2^T)(B^T-\hat B^T)\|\leq \epsilon_2 \|B^T\|\]

Then the final approximation is obtained as follows:
\[\tilde B^T=\hat B^T+ U_2U_2^T(B^T-\hat B^T)\]
\end{enumerate}
Now we verify that the approximation $\tilde B^T$ satisfies the constraints (\ref{eq:o2}), and (\ref{eq:o3}). 

First, we check whether $\hat B^T$ satisfies (\ref{eq:o2}).
\begin{align*}
(B^T-\hat B^T)X_w& =(B^T-U_1^s(U_1^s)^TB^TU_1^w(U_1^w)^T)X_w\\
&=(I-U_1^s(U_1^s)^T)B^TX_w+ U_1^s(U_1^s)^TB^T(I-U_1^w(U_1^w)^T)X_w
\end{align*}
Therefore, 
\[
\|(B^T-\hat B^T)X_w\|\leq  \epsilon_1(\|[X_s,B^TX_w]\|+ \|B^T\|\|[X_w,BX_s]\|)
\]
We can get a similar estimate for $X_s$:
\[
\|(B-\hat B)X_s\|\leq \epsilon_1(\|B^T\|\|[X_s,B^TX_w]\|+ \|[X_w,BX_s]\|)
\]

Next, we check $\tilde B^T$. Note that (\ref{eq:o3}) can be easily verified:
\begin{align*}
\|B^T-\tilde B^T\|&=\|(I-U_2U_2^T)(B^T-\hat B^T)\|\leq \epsilon_2\|B^T\|\\
\end{align*}
Now we check (\ref{eq:o2}) for $\tilde B^T$: 
\[
\|(B^T-\tilde B^T)X_w\| = \|(I-U_2U_2^T)(B^T-\hat B^T)X_w\|
\leq \epsilon_2\|B^T\| \|X_w\|
\]
We can also estimate it using $\epsilon_1$:
\begin{align*}
\|(B^T-\tilde B^T)X_w\|&= \|(I-U_2U_2^T)(B^T-\hat B^T)X_w\|\\
&\leq\|(B^T-\hat B^T)X_w\|\\
&\leq \epsilon_1(\|[X_s,B^TX_w]\|+ \|B^T\|\|[X_w,BX_s]\|)\\
\end{align*}

Therefore, we can only conclude that the upper bound for $\|(B^T-\tilde B^T)X_w\|$ satisfies: 
\[
\|(B^T-\tilde B^T)X_w\|\leq\min\{ \epsilon_1(\|[X_s,B^TX_w]\|+ \|B^T\|\|[X_w,BX_s]\|),  \epsilon_2\|B^T\| \|X_w\|\}
\]

Similarly,
\[\|(B-\tilde B)X_s\|\leq\min\{ \epsilon_1(\|B^T\|\|[X_s,B^TX_w]\|+ \|[X_w,BX_s]\|),  \epsilon_2\|B^T\| \|X_s\|\}\]

\subsubsection*{Second-order symmetric projection}

There is another symmetric approach to satisfy (\ref{eq:o2}) and (\ref{eq:o3}). Perform two truncated SVDs:
\[X_s=U_1^s\Sigma_1^s (V_1^s)^T+\hat U_1^s\hat\Sigma_1^s(\hat V_1^s)^T 
\qquad X_w=U_1^w\Sigma_1^w (V_1^w)^T+\hat U_1^w\hat\Sigma_1^w(\hat V_1^w)^T
\]
which satisfy 
\[\|\hat U_1^s\hat\Sigma_1^s(\hat V_1^s)^T\|\leq \epsilon_1\|X_s\|
\qquad 
\|\hat U_1^w\hat\Sigma_1^w(\hat V_1^w)^T\|\leq \epsilon_1\|X_w\|\]

Note that equivalently we have
\begin{equation}
\label{eq:error_sw}
\|(I-U_1^s(U_1^s)^T)X_s\| \leq \epsilon_1\|X_s\|
\qquad 
\|(I-U_1^w(U_1^w)^T)X_w\|\leq \epsilon_1\|X_w\|
\end{equation}

Then an intermediate approximation is \[\hat B^T=B^T- (I-U_1^s(U_1^s)^T)B^T(I-U_1^w(U_1^w)^T)\]

Perform another SVD:
\[B^T-\hat B^T= U_2\Sigma_2 V_2^T+E_2\]
where $\|E_2\|\leq \epsilon_2 \|B^T\|.$  

Then the final approximation is obtained as follows:
\[\tilde B^T=\hat B^T+ U_2U_2^T(B^T-\hat B^T)\]

Now we verify that the approximation $\tilde B^T$ satisfies the objectives (\ref{eq:o2}), and (\ref{eq:o3}). 

Note that
\begin{align*}
B^T-\tilde B^T&=(I-U_2U_2^T)(B^T-\hat B^T)\\
&=(I-U_2U_2^T)(I-U_1^s(U_1^s)^T)B^T(I-U_1^w(U_1^w)^T)
\end{align*}
We start with (\ref{eq:o3}): 
\begin{equation}
\label{eq:error_B}
\|B^T-\tilde B^T\|=\|E_2\|\leq \epsilon_2\|B^T\|
\end{equation}

Then we verify (\ref{eq:o2}).  For the ``s'' component, we have for $v_s=X_s\alpha$
\begin{align*}
\|(B-\tilde B) X_s\|&=\|(I-U_1^w(U_1^w)^T)B(I-U_1^s(U_1^s)^T)(I-U_2U_2^T)X_s\|\\
&=\|(I-U_1^w(U_1^w)^T)B(I-U_1^s(U_1^s)^T)(I-U_2U_2^T)(I-U_1^s(U_1^s)^T)X_s\|\\
&\leq \|B-\tilde B\| \|(I-U_1^s(U_1^s)^T)X_s\|\\
&\leq \epsilon_1\epsilon_2\|B^T\|  \|X_s\|
\end{align*}
Note that we always have
\[(I-U_1^s(U_1^s)^T)(I-U_2U_2^T)= (I-U_1^s(U_1^s)^T)(I-U_2U_2^T)(I-U_1^s(U_1^s)^T)\]
no matter whether $U_2$ and $U_1^s$ are orthogonal or not.

In particular, for one of the low-frequency eigenvectors $[e_i]_s$, $1\leq i\leq N$, we have
\[[e_i]_s=X_s\alpha_i
\qquad \text{where} \qquad
\alpha_i=
\begin{pmatrix}
0\\
\vdots\\
\lambda_i\\
\vdots\\
0
\end{pmatrix}
\]
Then \[\|(B-\tilde B)[e_i]_s\|\leq \lambda_i\epsilon_1\epsilon_2\|B^T\|\|X_s\|,\quad 1\leq i\leq N\]

For ``w'', it is similar.  We take $v_w=X_w\alpha$
\begin{align*}
\|(B^T-\tilde B^T) X_w\| &= \|(I-U_2U_2^T)(I-U_1^s(U_1^s)^T)B^T(I-U_1^w(U_1^w)^T)X_w\|\\
&=\|(I-U_2U_2^T)(I-U_1^s(U_1^s)^T)B^T(I-U_1^w(U_1^w)^T)(I-U_1^w(U_1^w)^T)X_w\|\\
&\leq \|B^T-\tilde B^T\| \|(I-U_1^w(U_1^w)^T)X_w\|\\
&\leq \epsilon_1\epsilon_2\|B^T\|  \|X_w\|
\end{align*}
Then 
\[\|(B^T-\tilde B^T)[e_i]_w\|\leq \lambda_i\epsilon_1\epsilon_2\|B^T\|\|X_w\|,\quad 1\leq i\leq N.\]

All these approaches approximately preserve vectors for hierarchical solvers.
\end{appendix}

\end{document}